\documentclass[12pt]{amsart}
\usepackage{amscd,amsmath,amssymb,amsfonts}
\usepackage[cmtip, all]{xy}

\newtheorem{thm}{Theorem}[section]
\newtheorem{prop}[thm]{Proposition}
\newtheorem{lem}[thm]{Lemma}
\newtheorem{cor}[thm]{Corollary}

\theoremstyle{definition}

\newtheorem{defn}[thm]{Definition}
\theoremstyle{remark}
\newtheorem{remk}[thm]{Remark}
\newtheorem{remks}[thm]{Remarks}

\newtheorem{exm}[thm]{Example}
\newtheorem{exms}[thm]{Examples}
\newtheorem{notat}[thm]{Notation}
\numberwithin{equation}{section}

{\hfill$\square$\end{defn}}

{\hfill$\square$\end{remk}}

{\hfill$\square$\end{remks}}

{\hfill$\square$\end{exm}}

{\hfill$\square$\end{exms}}

{\hfill$\square$\end{notat}}

\newcommand{\sC}{{\mathcal C}}
\newcommand{\sD}{{\mathcal D}}

\newcommand{\sH}{{\mathcal H}}

\newcommand{\sM}{{\mathcal M}}

\newcommand{\sQ}{{\mathcal Q}}
\newcommand{\sR}{{\mathcal R}}
\newcommand{\sS}{{\mathcal S}}
\newcommand{\sT}{{\mathcal T}}
\newcommand{\sU}{{\mathcal U}}
\newcommand{\sV}{{\mathcal V}}
\newcommand{\sW}{{\mathcal W}}
\newcommand{\sX}{{\mathcal X}}
\newcommand{\sY}{{\mathcal Y}}
\newcommand{\sZ}{{\mathcal Z}}
\newcommand{\A}{{\mathbb A}}
\newcommand{\B}{{\mathbb B}}
\newcommand{\C}{{\mathbb C}}
\newcommand{\D}{{\mathbb D}}
\newcommand{\E}{{\mathbb E}}
\newcommand{\F}{{\mathbb F}}
\newcommand{\G}{{\mathbb G}}
\renewcommand{\H}{{\mathbb H}}
\newcommand{\I}{{\mathbb I}}
\newcommand{\J}{{\mathbb J}}
\newcommand{\M}{{\mathbb M}}
\newcommand{\N}{{\mathbb N}}
\renewcommand{\P}{{\mathbb P}}
\newcommand{\Q}{{\mathbb Q}}
\newcommand{\R}{{\mathbb R}}
\newcommand{\T}{{\mathbb T}}
\newcommand{\U}{{\mathbb U}}
\newcommand{\V}{{\mathbb V}}
\newcommand{\W}{{\mathbb W}}
\newcommand{\X}{{\mathbb X}}
\newcommand{\Y}{{\mathbb Y}}
\newcommand{\Z}{{\mathbb Z}}
\renewcommand{\1}{{\mathbb{S}}}

\newcommand{\fm}{{\mathfrak m}}
\newcommand{\fM}{{\mathfrak M}}
\newcommand{\an}{{\rm an}}
\newcommand{\alg}{{\rm alg}}
\newcommand{\cl}{{\rm cl}}
\newcommand{\Alb}{{\rm Alb}}
\newcommand{\CH}{{\rm CH}}
\newcommand{\mc}{\mathcal}
\newcommand{\mb}{\mathbb}
\newcommand{\surj}{\twoheadrightarrow}
\newcommand{\inj}{\hookrightarrow}
\newcommand{\red}{{\rm red}}
\newcommand{\codim}{{\rm codim}}
\newcommand{\rank}{{\rm rank}}
\newcommand{\Pic}{{\rm Pic}}
\newcommand{\Div}{{\rm Div}}
\newcommand{\Hom}{{\rm Hom}}
\newcommand{\im}{{\rm im}}
\newcommand{\Spec}{{\rm Spec \,}}
\newcommand{\sing}{{\rm sing}}
\newcommand{\Char}{{\rm char}}
\newcommand{\Tr}{{\rm Tr}}
\newcommand{\Gal}{{\rm Gal}}
\newcommand{\Min}{{\rm Min \ }}
\newcommand{\Max}{{\rm Max \ }}
\newcommand{\supp}{{\rm supp}\,}
\newcommand{\0}{\emptyset}
\newcommand{\sHom}{{\mathcal{H}{om}}}

\newcommand{\Nm}{{\operatorname{Nm}}}
\newcommand{\NS}{{\operatorname{NS}}}
\newcommand{\id}{{\operatorname{id}}}
\newcommand{\Zar}{{\text{\rm Zar}}} 
\newcommand{\Ord}{{\mathbf{Ord}}}
\newcommand{\FSimp}{{{\mathcal FS}}}
\newcommand{\Sch}{{\operatorname{\mathbf{Sch}}}} 
\newcommand{\cosk}{{\operatorname{\rm cosk}}} 
\newcommand{\sk}{{\operatorname{\rm sk}}} 
\newcommand{\subv}{{\operatorname{\rm sub}}}
\newcommand{\bary}{{\operatorname{\rm bary}}}
\newcommand{\Comp}{{\mathbf{SC}}}
\newcommand{\IComp}{{\mathbf{sSC}}}
\newcommand{\Top}{{\mathbf{Top}}}
\newcommand{\holim}{\mathop{{\rm holim}}}
\newcommand{\Holim}{\mathop{{\it holim}}}
\newcommand{\op}{{\text{\rm op}}}
\newcommand{\<}{\langle}
\renewcommand{\>}{\rangle}
\newcommand{\Sets}{{\mathbf{Sets}}}
\newcommand{\del}{\partial}
\newcommand{\fib}{{\operatorname{\rm fib}}}
\renewcommand{\max}{{\operatorname{\rm max}}}
\newcommand{\bad}{{\operatorname{\rm bad}}}
\newcommand{\Spt}{{\mathbf{Spt}}}
\newcommand{\Spc}{{\mathbf{Spc}}}
\newcommand{\Sm}{{\mathbf{Sm}}}
\newcommand{\SmProj}{{\mathbf{SmProj}}}
\newcommand{\cofib}{{\operatorname{\rm cofib}}}
\newcommand{\hocolim}{\mathop{{\rm hocolim}}}
\newcommand{\Glu}{{\mathbf{Glu}}}
\newcommand{\can}{{\operatorname{\rm can}}}
\newcommand{\Ho}{{\mathbf{Ho}}}
\newcommand{\GL}{{\operatorname{\rm GL}}}
\newcommand{\SL}{{\operatorname{\rm SL}}}
\newcommand{\sq}{\square}
\newcommand{\Ab}{{\mathbf{Ab}}}
\newcommand{\Tot}{{\operatorname{\rm Tot}}}
\newcommand{\loc}{{\operatorname{\rm s.l.}}}
\newcommand{\HZ}{{\operatorname{\sH \Z}}}
\newcommand{\Cyc}{{\operatorname{\rm Cyc}}}
\newcommand{\cyc}{{\operatorname{\rm cyc}}}
\newcommand{\RCyc}{{\operatorname{\widetilde{\rm Cyc}}}}
\newcommand{\Rcyc}{{\operatorname{\widetilde{\rm cyc}}}}
\newcommand{\Sym}{{\operatorname{\rm Sym}}}
\newcommand{\fin}{{\operatorname{\rm fin}}}
\newcommand{\SH}{{\operatorname{\sS\sH}}}
\newcommand{\Wedge}{{\Lambda}}
\newcommand{\eff}{{\operatorname{\rm eff}}}
\newcommand{\rcyc}{{\operatorname{\rm rev}}}
\newcommand{\DM}{{\operatorname{\mathcal{DM}}}}
\newcommand{\GW}{{\operatorname{\rm{GW}}}}
\newcommand{\End}{{\operatorname{\text{End}}}}
\newcommand{\sSets}{{\mathbf{Spc}}}
\newcommand{\Nis}{{\operatorname{Nis}}}
\newcommand{\et}{{\text{\'et}}}
\newcommand{\Cat}{{\mathbf{Cat}}}
\newcommand{\ds}{{/\kern-3pt/}}
\newcommand{\bD}{{\mathbf{D}}}
\newcommand{\bC}{{\mathbf{C}}}
\newcommand{\res}{{\operatorname{res}}}

\renewcommand{\log}{{\operatorname{log}}}

\newcommand{\qf}{{\operatorname{q.fin.}}}

\newcommand{\fil}{\phi}
\newcommand{\barfil}{\sigma}

\newcommand{\Fac}{{\mathop{\rm{Fac}}}}
\newcommand{\Fun}{{\mathbf{Func}}}

\newcommand{\cp}{\coprod}
\newcommand{\ess}{\text{\rm{ess}}}
\newcommand{\dimrel}{\operatorname{\text{dim-rel}}}

 \newcommand{\gen}{{\text{gen}}}
\newcommand{\Gr}{{\text{\rm Gr}}}
\newcommand{\Ind}{{\operatorname{Ind}}}
\newcommand{\Supp}{{\operatorname{Supp}}}
\newcommand{\Cone}{{\operatorname{Cone}}}
\newcommand{\Tor}{{\operatorname{Tor}}}
\newcommand{\cok}{{\operatorname{coker}}}
\newcommand{\Br}{{\operatorname{Br}}}
\newcommand{\sm}{{\operatorname{sm}}}
\newcommand{\Proj}{{\operatorname{Proj}}}
\newcommand{\Pos}{{\operatorname{\mathbf{Pos}}}}
\newcommand{\colim}{\mathop{\text{colim}}}
\newcommand{\ps}{{\operatorname{\text{p.s.s.}}}}
\newcommand{\Pro}{{\mathbf{Pro}\hbox{\bf -}}}
\newcommand{\hocofib}{{\operatorname{\rm hocofib}}}
\newcommand{\hofib}{{\operatorname{\rm hofib}}}
\newcommand{\Th}{{\operatorname{\rm Th}}}
\newcommand{\Cube}{{\mathbf{Cube}}}
\newcommand{\Mod}{{\mathbf{Mod}}}
\newcommand{\alt}{{\operatorname{\rm alt}}}
\newcommand{\Alt}{{\operatorname{\rm Alt}}}
\newcommand{\sgn}{{\operatorname{\rm sgn}}}
\newcommand{\dga}{{\rm d.g.a.\!}}
\newcommand{\cdga}{{\rm c.d.g.a.\!}}
\newcommand{\lci}{{\rm l.c.i.\!}}
\newcommand{\sym}{{\operatorname{\rm sym}}}
\newcommand{\DTM}{{\operatorname{DTM}}}
\newcommand{\DCM}{{\operatorname{DCM}}}
\newcommand{\MTM}{{\operatorname{MTM}}}
\newcommand{\Mot}{{\operatorname{Mot}}}
\renewcommand{\sp}{{\widetilde{sp}}}
\newcommand{\NST}{{\mathbf{NST}}}
\newcommand{\PST}{{\mathbf{PST}}}
\newcommand{\Sus}{{\operatorname{Sus}}}
\newcommand{\tr}{{\operatorname{tr}}}
\newcommand{\cone}{{\operatorname{cone}}}
\newcommand{\na}{{\rm{naive}}}
\newcommand{\Cor}{{\operatorname{Cor}}}
\newcommand{\DMH}{\operatorname{DM}^H}
\newcommand{\gm}{\text{gm}}
\newcommand{\TZ}{{\operatorname{Tz}}}
\renewcommand{\TH}{{\operatorname{TCH}}}
\newcommand{\un}{\underline}
\newcommand{\ov}{\overline}
\renewcommand{\hom}{\text{hom}}
\newcommand{\dgn}{{\operatorname{degn}}}
\renewcommand{\dim}{\text{\rm dim}}
\newcommand{\td}{\mathbf{d}}
\newcommand{\tuborg}{\left\{\begin{array}{ll}}
\newcommand{\sluttuborg}{\end{array}\right.}
\newcommand{\dto}{\dashrightarrow}

\begin{document}
\title{Mixed motives over $k[t]/{(t^{m+1})}$}
\author{Amalendu Krishna, Jinhyun Park}
\address{School of Mathematics, Tata Institute of Fundamental Research,  
1 Homi Bhabha Road, Colaba, Mumbai, India}
\email{amal@math.tifr.res.in}
\address{Department of Mathematical Sciences, KAIST,  
Daejeon, 305-701, Republic of Korea (South)}
\email{jinhyun@mathsci.kaist.ac.kr; jinhyun@kaist.edu}

\baselineskip=10pt 
  

\maketitle
\begin{abstract}For a perfect field $k$, we use the techniques of 
Bondal-Kapranov 
and Hanamura to construct a triangulated category of mixed motives 
over the truncated polynomial ring $k[t]/(t^{m+1})$. The extension groups
in this category are given by Bloch's higher Chow groups and the additive
higher Chow groups. The main new ingredient is the moving lemma 
for additive higher Chow groups in and its refinements.
\end{abstract}

\section{Introduction}
Let $k$ be perfect field. The aim of this paper is to construct a 
triangulated category of mixed motives over the truncated polynomial ring 
$k[t]/(t^{m+1})$, such that the resulting motivic cohomology groups 
(the Ext groups) of smooth projective varieties in this category 
compute the $K$-theory of perfect complexes on the infinitesimal 
deformations of these varieties. In other words, this category is expected
to be an extension of the category of mixed motives over $k$,
constructed for example in \cite{Ha1}, \cite{Levine} and \cite{Voevodsky}, to
the simplest types of non-reduced rings. The complete construction of such
a category with expected properties has been desired for a long time 
(\emph{cf.} \cite{Bloch-tangent}) and
will go a long way in understanding how one could construct 
the motivic cohomology that compute the $K$-theory of vector bundles
on singular varieties. In this paper, our focus is to study such a problem
in the particular case of those singular varieties which are the 
infinitesimal deformations of smooth varieties.
  
Following Bloch's proposal on how to construct mixed Tate motives over the 
field $k$ in \cite{Bloch}, it was conjectured by Bloch and Esnault in 
\cite{BE2} that there should be a theory of \emph{``degenerate''} cycle 
complexes over $k$ in such a way that the Tanakian formalism of \cite{Bloch}
could be used to construct the category of mixed Tate motives over the
ring of dual numbers $k_{\epsilon} := k[t]/{(t^2)}$. 

Possibly motivated by 
\cite{BE2}, Goncharov \cite{Goncharov} used his idea of $k$-scissors 
congruence to define Euclidean scissors congruence groups to get a
Lie coalgebra ${\sQ}_{\bullet}(k_{\epsilon})$ in the category of 
${\Q}_{\epsilon}$-modules. He conjectured that if $k$ is algebraically closed,
the category of finite dimensional graded comodules over 
${\sQ}_{\bullet}(k_{\epsilon})$ should be equivalent to a subcategory
of the conjectured category of mixed Tate motives over $k_{\epsilon}$
(\emph{cf.} \cite[Conjecture~1.3]{Goncharov}). However, a complete
construction of even the category of these mixed Tate motives was not
known so far.

Our aim in this paper is to give a complete construction
of a bigger category of mixed motives over any given truncated
polynomial ring $k_m = k[t]/(t^{m+1})$ (note that $k_0 = k$).
We expect that this category has
all the expected properties. In particular, the Ext groups should give the
$K$-theory of the infinitesimal thickenings of smooth projective varieties.
Although we are unable to prove this last property, there are strong 
indications that this should indeed be true as we shall see shortly. Some more 
results in this direction will appear in \cite{KP1}. 

Before we describe the main result,
we fix some terminology. 
Let $\SmProj/k$ denote the category of smooth projective varieties over $k$.
The category of all quasi-projective schemes over $k$ will be denoted 
by $\Sch/k$. The subcategory with only proper
morphisms will be denoted by ${\Sch}'/k$. Let $\mathcal{DM}(k)$ denote 
(the integral version of)
Hanamura's triangulated category of mixed motives over $k$ ({\sl cf.} 
\cite{KL}). For a 
$X \in \Sch/k$, let $\TH^r_{\log}(X, n;m)$ denote the log
additive higher Chow groups as in \cite{KP} (see also \cite{KL}). Note
that for $X$ smooth and projective, these are just the additive higher 
Chow groups $\TH^r (X, n; m)$. We refer to 
Section~\ref{section:cyclecomplex} for the
review of these groups. Let $\CH^r(X, n)$ denote the higher Chow groups
of $X$. Finally, recall from \cite{Soule} that for a scheme $X$, the higher
$K_i(X)$-groups of perfect complexes on $X$ have gamma filtrations which
induces Adams operations on each ${K_i(X)}_{\Q}$. The $r$-th eigenspace for
this operation is denoted by $K^{(r)}_i(X)$. 
\begin{thm}\label{thm:Main}
For $m \ge 0$, there exists a triangulated category 
$\mathcal{DM}(k;m)$ such that the following results hold:
\begin{enumerate}
\item  There are natural functors
$$\iota: \mathcal{DM}(k) \to \mathcal{DM}(k;m),$$ which is faithful (but not 
full) and 
$${\rm Forget}: \mathcal{DM}(k;m) \to \mathcal{DM}(k)$$ 
such that ${\rm Forget} \circ \iota$ is the identity. Moreover, 
$\mathcal{DM}(k;0)$ is canonically isomorphic to $\mathcal{DM}(k)$. 
\item There exists the motive functor with the modulus $m$ augmentation, 
$$h : \SmProj/k \to \mathcal{DM}(k;m)$$ such that 
$$\Hom_{\mathcal{DM}(k; m)} (\un{\mathbb{Z}}, h(X)(r)[2r-n]) = 
\CH^r (X, n) \oplus \TH^r (X, n;m),$$
where $\un{\mathbb{Z}}(r)$ and $(-) (r)$ are Tate objects and Tate twists.
\item If $k$ admits Hironaka's resolution of singularities, 
then there exists an extension of $h(-)$
$$bm : {\Sch}'/k \to \mathcal{DM}(k;m)$$ 
such that
\[
\Hom_{\mathcal{DM}(k; m)} (\un{\mathbb{Z}}, bm(X)(r)[2r-n]) = 
\CH^r (X, n) \oplus \TH^r_{\log}(X, n;m)
\]
for a smooth quasi-projective variety $X$.
\end{enumerate}
\end{thm}
As the reader will observe, the motive functor $h$ in the above theorem is
a covariant functor unlike the one in \cite{Ha1}, which is a contravariant
functor.
By combining the results of Hesselholt \cite{Hesselholt} and 
\cite[Theorems~3.4, 3.7]{KP} (see also \cite{R}), we have the following
consequence of the above theorem:
\begin{cor}\label{cor:morph*}
$1$. For $n \ge 1$, there is an isomorphism 
\[
\Hom_{\mathcal{DM}(k; m)} (\un{\mathbb{Q}}, h({\rm Spec}(k))(n)[n]) \cong 
K^{(n)}_n\left(k[t]/{(t^{m+1})}\right).
\]
$2$. Suppose $k$ is algebraically closed of characteristic zero. For 
$n \ge 3$, there is a natural surjection 
\[ 
\Hom_{\mathcal{DM}(k; 1)} (\un{\mathbb{Q}}, h({\rm Spec}(k))(n-1)[n-2]) \surj 
K^{(n-1)}_n\left(k[t]/{(t^{2})}\right).
\]
\end{cor}
As seen in \cite{KL}, the additive higher Chow groups are expected
to give rise to Atiyah-Hirzebruch spectral sequence
\[
\TH^{-q}_{\rm log}(X, -p-q;m) \Rightarrow K^{\rm log}_{-p-q}(X),
\]
where $K^{\rm log}$ is a spectrum which for $X \in \SmProj/k$, is the
homotopy fiber of the map of spectra $K\left(X[t]/{(t^{m+1})}\right) \to K(X)$.
Based on this and  the above computations and the ongoing work \cite{KP1},
we expect that for $X \in \SmProj/k$ and for $n \ge 1$, there is a natural 
isomorphism
\begin{equation}\label{eqn:Ktheory}
\Hom_{\mathcal{DM}(k; m)} (\un{\mathbb{Q}}, h(X)(r)[2r-n]) = 
{K^{(r)}_n\left(X[t]/{(t^{m+1})}\right)}.
\end{equation}
 
We also expect Goncharov's category \cite{Goncharov} to be a subcategory
of $\mathcal{MTM}(k; m)$ ({\sl cf.} ~\eqref{subsection:Ttwist}). We hope
that this can be proved using the techniques in \cite{Bloch} and \cite{KP},
where we showed that the additive higher Chow groups have a natural
structure of differential graded algebra. In a sequel to this work, we shall
address the question of generalizing Voevodsky's category of mixed motives 
over $k$ to the truncated polynomial ring $k[t]/{(t^{m+1})}$.

We now give a brief outline of this paper. Our construction of the
category of motives is broadly based on the construction of triangulated
categories out of a dg-category in Bondal-Kapranov \cite{BK} and the
construction of ${\sD}{\sM}(k)$ in Hanamura \cite{Ha1}. In order to carry
this out, we formalize the results of \cite{BK} and \cite{Ha1} in the
language of what we call a \emph{partial} dg-category, in the next two 
sections. 
Apart from its use in this
paper, we hope that this formalism of partial dg-categories will be useful
in proving many other similar results, especially in constructing various
types of categories of motives. We review (additive) higher Chow groups and 
some of their properties in Section~\ref{section:cyclecomplex}.
Section~\ref{section:ML} contains the main technical results about the
moving lemma and its refinements for additive higher Chow groups. 
We construct our category $\mathcal{DM}(k;m)$ in Section~\ref{section:DM} 
using the results of 
Section~\ref{section:TPDG} and ~\ref{section:ML}. The last section contains 
results on the
extension of the motives to all schemes of finite type over $k$. This is
based on some results of \cite{GN} and \cite{KL}.

\section{Partial dg-category}\label{section:PDG}
The construction of the triangulated category ${\sD}{\sM}(k;m)$ of mixed
motives over $k[t]/{(t^{m+1})}$ in this paper is broadly based
on a very general construction of a triangulated category from a dg-category,
by Bondal and Kapranov in \cite{BK}. Apart from \emph{ibid.}, our construction
crucially uses the modification of Bondal-Kapranov's techniques by
Hanamura in the construction of his category of mixed motives in \cite{Ha1}.

Given a \emph{pre-additive} dg-category $\sC$, Bondal and Kapranov construct
a sequence of dg-categories and functors
\[
\sC \to {\sC}^{\oplus} \to {\rm {PreTr}}(\sC) \to {\rm Tr}(\sC),
\]
such that all intermediate categories are additive dg-categories and the end 
product ${\rm Tr}(\sC)$ is a triangulated category. Moreover, this
construction is natural with respect to functors of pre-additive
dg-categories.

It turns out that this formalism of Bondal-Kapranov can be adapted
also to slightly more general setting where one allows more flexibility
on the composability axioms about the morphisms in the underlying
dg-categories. It is this refinement of the construction of \cite{BK}
that will be needed to obtain the category ${\sD}{\sM}(k;m)$.

In this section, we carry out the construction of Bondal and Kapranov
in this more general setting of what we shall call 
\emph{partial dg-categories}. This
new formalism of partial dg-categories is motivated by the construction
of mixed motives in \cite{Ha1}. In fact, our endeavor in this and the next
section is to axiomatize the techniques of \emph{ibid.} in the language of
partial dg-categories. We hope that this abstract formalization will
be useful in many cases of interest, especially where one works with 
algebraic cycles and motives.

Let $K(\Z)$ and $D(\Z)$ respectively denote the homotopy category of cochain 
complexes of abelian groups and its derived category. Similarly, let $K^{-}(\Z)$ and $D^{-}(\Z)$ denote the corresponding categories of right 
bounded cochain complexes. 
We shall say that $f: M^{\bullet} \dto N^{\bullet}$ is a 
\emph{partially defined} morphism of cochain complexes in $K(\Z)$, if there is a
subcomplex ${M'}^{\bullet} \overset{i}{\inj} {M}^{\bullet}$, where
$i$ is a quasi-isomorphism, and $f:{M'}^{\bullet} \to N^{\bullet}$ is an 
honest morphism of cochain complexes.  For a cochain complex $M^{\bullet}$, 
the term {\sl quasi-isomorphic subcomplex} will mean a subcomplex  
${M'}^{\bullet} \inj M^{\bullet}$ such that the inclusion is a 
quasi-isomorphism.  

Recall that a dg-category over $\Z$ consists of an
additive category $\sT$ such that 

\noindent $1.$ For any pair of objects $A, B$ in $\sT$, one has $\left(\Hom_{\sT}(A,B),
d\right) \in K(\Z)$. 

\noindent $2.$ For a triple $(A,B,C)$ of objects in $\sT$, there is a composition
morphism ${\mu}_{ABC}: \Hom_{\sT}(A,B) {\otimes}_{\Z} \Hom_{\sT}(B,C)
\to \Hom_{\sT}(A,C)$.

\noindent $3.$ For any object $A$ of $\sT$, there is a ``unit'' morphism  
$e_A : \Z \to \Hom_{\sT}(A,A)$ of complexes. \\
The composition of morphisms satisfies the usual associative law and the
left and right compositions with the unit morphism $e_A$ act as identity on 
any $\Hom_{\sT}(A,B)$ and $\Hom_{\sT}(B,A)$. A dg-category $\sT$ as above
which does not necessarily have finite coproducts of objects is called
a pre-additive dg-category.  

In what follows, we consider a more general analogue of a dg-category,    
where the compositions of morphisms are only partially defined in the
above sense.

\begin{defn}\label{defn:Pdg}
A \emph{partial dg-category} $\sC$ over $\Z$ consists of the following data. \\
$({\rm P}1)$ A set of objects $Ob(\sC)$, also denoted by $\sC$ itself. \\ 
$({\rm P}2)$ For any pair of objects $A, B$ in $\sC$, one has $\Hom_{\sC}(A,B)
\in K(\Z)$, and a collection $S(A,B)$ of quasi-isomorphic subcomplexes of
$\Hom_{\sC}(A,B)$ called
``\emph{distinguished subcomplexes}''. \\
$({\rm P}3)$ For any object $A$ of $\sC$, there is a ``unit'' morphism  
$e_A : \Z \to \Hom_{\sC}(A,A)$ of complexes. \\
$({\rm P}4)$ Given any 
$f \in \Hom_{\sC}(A,B), g \in \Hom_{\sC}(B,C)$ and a
distinguished subcomplex
\[
{\Hom_{\sC}(A,C)}' \subset \Hom_{\sC}(A,C),
\]
there are distinguished subcomplexes ${\Hom_{\sC}(B,C)}' \subset 
\Hom_{\sC}(B,C)$, ${\Hom_{\sC}(A,B)}' \subset \Hom_{\sC}(A,B)$ such that the 
compositions 
\[
 (-) \circ f:{\Hom_{\sC}(B,C)}'  \to  {\Hom_{\sC}(A,C)}', \ {\rm and}
\]
\[
g \circ (-): {\Hom_{\sC}(A,B)}' \to {\Hom_{\sC}(A,C)}'
\]
are defined. \\  
$({\rm P}5)$ For any pair of objects $A, B$ in $\sC$ and for any two
distinguished subcomplexes $M, M' \subset \Hom_{\sC}(A,B)$, there is a
distinguished subcomplex $M'' \subset M \cap M'$ of $ \Hom_{\sC}(A,B)$. \\
$({\rm P}6)$ The composition of morphisms at the level of distinguished
subcomplexes satisfies the associative law, and the
partially defined left and right compositions with the unit morphism $e_A$ 
act as identity on any $\Hom_{\sC}(A,B)$ and $\Hom_{\sC}(B,A)$.
\end{defn}
A partial dg-category $\sC$ which has all finite coproducts of its objects 
will be called an \emph{additive} partial dg-category.
An example of a partial dg-category will be given later in this paper
when we construct our category ${\sD}{\sM}(k;m)$. 

If $\sC$ is a partial dg-category, let ${\sC}^{\oplus}$ be the
partial dg-category whose objects are formal finite coproducts of the 
objects of $\sC$, i.e., 
\[
A = \bigoplus_{u \in J} A_u
\]
where $A_{u} \in Ob(\sC)$ and $|J|<\infty$. If 
$J = \emptyset$, then we write $A=0 $ by convention. It is easy to see
that for the possibly pre-additive partial dg-category $\sC$, the new category
${\sC}^{\oplus}$ is indeed an additive partial dg-category. If $\mathcal{C}$ is additive from the first place, then $\mathcal{C} = \mathcal{C}^{\oplus}$.

\subsection{Twisted complexes and ${\rm PreTr}(\sC)$}
Let $\sC$ be a partial dg-category. 
\begin{defn}[{\cite{BK}}]\label{defn:tComplex}
A \emph{twisted complex over} $\mathcal{C}$ is a system  
$A = \{ (A^i)_{i \in \mathbb{Z}}, q_{i,j} : 
A^i \to A^{j} \ {\rm for} \ i < j\}$, where \\
$\bullet$ $A^i \in Ob(\mathcal{C}^{\oplus})$, all but finitely many of them 
are $0$, and $q_{i,j}$ are morphisms in $\mathcal{C}^{\oplus}$ of degree 
$i - j + 1$. \\
$\bullet$ For any sequence $i = i_0 < \cdots < i_r = j$, the compositions
$q_{i_{r-1}, i_{r}} \circ \cdots \circ q_{i_0, i_1}$ are defined. \\
$\bullet$ For all $i<j$,
\begin{equation}\label{eqn:tComplex1}
(-1)^j d (q_{i,j}) + \stackrel{}{\underset {i<k<j}{\sum}} q_{k,j} 
\circ q_{i,k} = 0.
\end{equation} 
\end{defn}
Note that the twisted complexes as defined above are analogous to the 
one-sided twisted complexes of \cite[Definition~4.1]{BK}. 
\begin{remk}\label{remk:finite1}
Note also that since only finitely many $A^i$'s are non-zero in a 
twisted complex $A$ and since there is exactly one given 
$q_{i,j} : A^i \to A^j$, the system
$A$ involves only finitely many nonzero morphisms $q_{i,j}$'s, too. In particular, if $A^i = 0$ for all but one $i$, then all $q_{i,j}=0$.
\end{remk}
We now define the set of partial morphisms between two twisted
complexes. So let $A = \{ (A^i)_{i \in \mathbb{Z}}, p_{i,j} : 
A^i \to A^{j} \}, B = \{ (B^{i'})_{i' \in \mathbb{Z}}, q_{i',j'} : 
B^{i'} \to B^{j'} \}$ be two twisted complexes over $\sC$.  Write
$A^i = \stackrel{}{\underset {\alpha \in I(i)}{\oplus}} A_{\alpha}^i$
and $B^{i'} = 
\stackrel{}{\underset {\beta \in I'(i')}{\oplus}} B_{\beta}^{i'}$.
The axiom $({\rm P}5)$ of the definition of a partial dg-category implies that
given any finite collection $M_i \subset \Hom_{\sC}(A,B)$ of
distinguished subcomplexes, there is a distinguished subcomplex
$M \subset \left(\stackrel{}{\underset {i}{\cap}} M_i\right)$ of 
$\Hom_{\sC}(A,B)$. Using this, the axiom $({\rm P}4)$ of Definition~\ref{defn:Pdg},
and Remark~\ref{remk:finite1}, we can find distinguished subcomplexes
${\Hom_{\sC}(A_{\alpha}^i,B_{\beta}^{i'})}' \subset 
\Hom_{\sC}(A_{\alpha}^i,B_{\beta}^{i'})$ so that for
\[
{\Hom_{{\sC}^{\oplus}}(A^i,B^{i'})}' := 
\stackrel{}{\underset {\alpha}{\bigoplus}}
\stackrel{}{\underset {\beta}{\bigoplus}} 
{\Hom_{\sC}(A_{\alpha}^i,B_{\beta}^{i'})}',
\]
the following holds: for any sequence $i = i_0 < \cdots < i_r = j$,  $ i'= i'_0 <
\cdots < i'_s = j'$, the composition 
\[
{\Hom_{{\sC}^{\oplus}}(A^{j},B^{i'})}' \to 
{\Hom_{{\sC}^{\oplus}}(A^i,B^{j'})},
\]
\[
u \mapsto {\underbrace{q_{i'_{s-1}, i'_s} \circ 
\cdots q_{i'_0, i'_1}}} \circ u \circ
{\underbrace{p_{i_{r-1}, i_r} \circ \cdots \circ p_{i_0, i_1}}}
\]
is defined. The axiom $({\rm P}6)$ of Definition~\ref{defn:Pdg} then implies
that these compositions are associative. In particular, the maps
\[
(-) \circ p_{i,j} : {\Hom_{{\sC}^{\oplus}}(A^j,B^{i'})}' \to
{\Hom_{{\sC}^{\oplus}}(A^i, B^{i'})}
\]
are defined and so are $q_{i',j'} \circ (-)$. One defines the complex
$\Hom_{{\rm PreTr}(\sC)}(A,B)$ as the cochain complex 
\begin{equation}\label{eqn:PreTr*}
{\Hom_{{\rm PreTr}(\sC)}(A,B)}^n =
\stackrel{}{\underset {-i+j+l = n}{\bigoplus}}
{\left({\Hom_{{\sC}^{\oplus}}(A^i,B^{j})}'\right)}^l.
\end{equation}
The differential $D$ of the complex $\Hom_{{\rm PreTr}(\sC)}(A,B)$ is 
given for $f \in {\left({\Hom_{{\sC}^{\oplus}}(A^i,B^{j})}'\right)}^l$
by the formula
\begin{equation}\label{eqn:D} 
D (f) := (-1)^j d (f) + \sum_m\left( (-1)^{j+m} q_{j, m}  
\circ f + (-1)^{l+j+m+1} f \circ p_{m, i}\right).
\end{equation}
One should note here that the various signs in the differential are completely different from the ones chosen in \cite{BK} and they conform
more to Hanamura's construction.
\begin{lem}\label{lem:D2=0}
The above $D$ satisfies $D \circ D = 0.$
\end{lem}

\begin{proof} For $f \in (\hom_{\mathcal{C}^{\oplus}} (A^i, B^{j})')^l$, in the formula \eqref{eqn:D}, 
we let $(A) := (-1)^j d(f)$, $(B): = \sum_m (-1)^{j+m} q _{j,m} \circ f$, and 
$(C):=  \sum_m (-1)^{l+j+m+1} f \circ p_{m,i}$ so that 
$D(f) = (A) + (B) + (C)$. We prove that $D^2 (f) = D (A) + D(B) + D(C) = 0$.

First we have \begin{eqnarray*}
D(A) &=&(-1)^j d (df) +  \sum_m \left( (-1)^m q_{j,m} \circ 
df + (-1)^{l+m} df \circ p_{m,i} \right) \\
&=&  \underset{(A1)}{\underbrace{\sum_m (-1)^m q_{j,m} \circ df}} + 
\underset{(A2)}{\underbrace{\sum_m (-1)^{l+m} df \circ p_{m,i}}}\\
\end{eqnarray*}

For $(B)$, a direct calculation shows that
\begin{eqnarray*}
D(B)&=& \underset{(B1)}{\underbrace{\sum_m (-1)^j d (q_{j,m} \circ f) }}+  
\underset{(B2)}{\underbrace{ \sum_{m,m'} (-1)^{j+m'} q_{m,m'} \circ 
q_{j, m} \circ f }} \\
&+& \underset{(B3)}{\underbrace{ \sum_{m, m'} (-1)^{l+m+m'} q_{j,m} \circ 
f \circ p_{m', i}}},
\end{eqnarray*}where the Leibniz rule for $(B1)$ shows that we have

\begin{eqnarray*}
(B1) = \underset{(B11)}{\underbrace{ \sum_m (-1)^j dq_{j,m} \circ f}} + 
\underset{(B12)}{\underbrace{ \sum_{m} (-1)^{1-m} q_{j,m} \circ df}}.
\end{eqnarray*}

Similarly, a direct calculation shows that
\begin{eqnarray*}
D(C) &=& \underset{(C1)}{\underbrace{ \sum_{m} (-1)^{l+m+1} 
d( f \circ p_{m,i})}} + \underset{(C2)}
{\underbrace{ \sum_{m,m'} (-1)^{l+m+m'+1} q_{j, m'} \circ f \circ 
p_{m,i} }} \\
&+& \underset{(C3)}{\underbrace{ \sum_{m,m'} (-1)^{m'-i+1} f \circ p_{m,i} 
\circ q_{m',m}}},
\end{eqnarray*}where the Leibniz rule for $(C1)$ shows that we have

\begin{eqnarray*}
(C1)= \underset{(C11)}{\underbrace{ \sum_m (-1)^{l + m + 1} df \circ 
p_{m,i} }} + \underset{(C12)}{\underbrace{ \sum_m (-1)^{m+1} f \circ 
dp_{m,i}}}.
\end{eqnarray*}

Now, one immediately notices that $(A1) + (B12) = 0$, $(A2) + (C11) = 0$, 
$(B3) + (C2) = 0$, and $(B2) + (B11) = (C3) + (C12) = 0$ by the condition 
\eqref{eqn:tComplex1}. Thus, $D^2 (f) = D(A) + D(B) + D(C) = 0$. This 
proves the lemma.
\end{proof} 
\begin{defn}\label{defn:PreTRC}
Let $\sC$ be a partial dg-category. We define ${\rm PreTr}(\sC)$ to be
the \emph{partial} dg-category whose objects are all the twisted complexes
over $\sC$, and whose ``morphisms'' are given by the cochain complexes defined
in ~\eqref{eqn:PreTr*}.
\end{defn}
Observe that ${\rm PreTr}(\sC)$ is not yet an honest category since the
morphisms between twisted complexes depend on the choices of distinguished
subcomplexes. We shall show however that these morphisms are well defined
up to quasi-isomorphisms. \\
A \emph{full subcategory} $\sD$ of ${\rm PreTr}(\sC)$
is a partial category such that \\
$\bullet$ $Ob(\sD) \subset Ob\left({\rm PreTr}(\sC)\right)$ and \\
$\bullet$ For $A, B \in Ob(\sD)$, $\Hom_{\sD}(A,B) =
\Hom_{{\rm PreTr}(\sC)}(A,B)$. \\
In particular, the morphisms of $\sD$ will be shown to be well
defined up to quasi-isomorphisms.  

\section{Homotopy category of a partial dg-category}\label{section:TPDG}
In this section, we complete the program of constructing an honest 
triangulated category ${\rm Tr}(\sC)$ from a given partial dg-category $\sC$, 
which will be called the \emph{homotopy category of $\sC$}. This is done with the help
of the notion of $C$-complexes (\emph{cf.} \cite[Section~3]{Ha1}).
We begin with a brief recall of this theory. We call it a \emph{left} $C$-complex here.\\
\textbf{Notational convention}: In this section, when we write sums over 
various indices, we emphasize the indices over which the sums are taken by 
putting underlines for them. For example, $\sum_{m < \un{k} < n}$ is taken
over $k$ such that $m < k<n$ with $m,n$ fixed, while 
$\sum_{m < \un{k} < \un {n}}$ is taken over all pairs of indices $k$ and $n$
with $m<k<n$ for a fixed $m$.

\begin{defn}\label{defn:C-complex}
A \emph{left $C$-complex} of abelian groups consists of \\
$(i)$ A sequence of cochain complexes $\left(A^{\bullet}_m, d_{A_m}\right)$
for $m \in \Z$ such that $A^{\bullet}_m = 0$ for all but finitely many $m$'s.
\\
$(ii)$ For $m < n$, there are maps of graded groups
\[
F_{m,n} : A^{\bullet}_m \to A^{\bullet}_n [m - n +1]
\]
subject to the condition
\begin{equation}\label{eqn:C-complex1}
F_{m,n} \circ (-1)^m d_{A_m} + (-1)^nd_{A_n} \circ F_{m,n} +
\stackrel{}{\underset {m < \un{l} < n}{\sum}} F_{l,n} \circ F_{m,l} = 0
\end{equation}
as a map $A^{\bullet}_m \to A^{\bullet + m - n +2}_n$.
\end{defn}
Given a left $C$-complex $\left(A^{\bullet}_m, d_{A_m}\right)$, one defines
its total complex $Tot(A) = \left({Tot(A)}^{\bullet}, \bf{d}^L\right)$
by 
\begin{equation}\label{eqn:C-complex10}
{Tot(A)}^{p} = \stackrel{}{\underset {m+ \un{i} = p}{\bigoplus}} A^{i}_m
\end{equation}
such that for $f \in A^{p-m}_m$, one has 
\begin{equation}\label{eqn:C-complex2}
{\bf{d}^L}(f) = \left((-1)^m d_{A_m}(f) + 
\stackrel{}{\underset {\underline{n} > m}{\sum}} 
F_{m,n}(f)\right) \in \ 
\stackrel{}{\underset {\underline{n} \ge m}{\bigoplus}} 
A^{p-n+1}_{n}.
\end{equation}
One checks using the condition ~\eqref{eqn:C-complex1} that ${\bf{d}^L}$ is 
indeed a differential.
\begin{remk}\label{remk:C-complex3}
As explained in \emph{loc. cit.}, a left $C$-complex is a generalization of
the notion of double complexes in that the maps $F_{m, m+1}$ are
chain maps such that $F_{m+1, m+2} \circ F_{m, m+1}$ are \emph{not} assumed to be 
zero, although they are zero in the homotopy category $K(\Z)$ via the 
homotopy $F_{m, m+2}$. In fact, the maps $F_{m, n}$ of higher lengths
$(m-n)$ give the null-homotopy for the composites of the similar maps
of smaller lengths. In particular, a left $C$-complex is a chain complex
of objects in the homotopy category $K(\Z)$ of chain complexes.  The
standard formalism of spectral sequences associated to a chain complex
of objects in $K(\Z)$ then implies that there is a convergent spectral sequence
\begin{equation}\label{eqn:SS1}
E^{p,q}_1  = H^q(A^{\bullet}_p) \Rightarrow H^{p+q}\left(Tot(A),
{\bf{d}^L}\right).
\end{equation} 
\end{remk}
Apart from the above left $C$-complexes, we shall also need the following
variant of these objects that we shall call \emph{right} $C$-complexes. 
We emphasize again that the above left $C$-complexes are exactly what are
simply called $C$-complexes in \cite{Ha1}.
\begin{defn}\label{defn:RCC}
A \emph{right $C$-complex} of abelian groups consists of \\
$(i)$ A sequence of cochain complexes $\left(A^{\bullet}_m, d_{A_m}\right)$
for $m \in \Z$ such that $A^{\bullet}_m = 0$ for all but finitely many $m$'s.
\\
$(ii)$ For $m < n$, there are maps of graded groups
\[
E_{m,n} : A^{\bullet}_m \to A^{\bullet}_n [m - n +1]
\]
subject to the condition
\begin{equation}\label{eqn:RCC1}
E_{m,n} \circ (-1)^m d_{A_m} + (-1)^nd_{A_n} \circ E_{m,n} +
\stackrel{}{\underset {m < \underline{l} < n}{\sum}} 
(-1)^{l+1} E_{l,n} \circ E_{m,l} = 0
\end{equation}
as a map $A^{\bullet}_m \to A^{\bullet + m - n +2}_n$.
\end{defn}

Given a right $C$-complex $\left(A^{\bullet}_m, d_{A_m}\right)$, one defines
its total complex $Tot(A) = \left({Tot(A)}^{\bullet}, \bf{d}^R\right)$
by 
\begin{equation}\label{eqn:CR-complex10}
{Tot(A)}^{p} = \stackrel{}{\underset {m+\un{i} = p}{\bigoplus}} A^{i}_m
\end{equation}
such that for $f \in A^{p-m}_m$, one has 
\begin{equation}\label{eqn:CR-complex2}
{\bf{d}^R}(f) = \left((-1)^m d_{A_m}(f) + 
\stackrel{}{\underset {\underline{n} > m}{\sum}} 
(-1)^{n+1} E_{m,n}(f)\right) \in \ 
\stackrel{}{\underset {\underline{n} \ge m}{\bigoplus}} 
A^{p-n+1}_{n}.
\end{equation} 
\begin{lem}\label{lem:RCCTD}
${\bf{d}^R}\circ {\bf{d}^R} = 0$. In other words, $\left(Tot(A), 
{\bf{d}^R}\right)$ is a cochain complex.
\end{lem}
\begin{proof}
By a direct calculation, for $f \in A_m ^{p-m}$, 
\begin{eqnarray*}
{(\bf{d}^R)}^2 (f) &=& {\bf{d^R}} \left( (-1)^m d_A (f) + 
\sum_{\un{n} > m } (-1)^{n+1} E_{m,n} (f) \right) \\
 &=& \underset{(A)}{\underbrace{{\bf{d^R}} \left((-1)^m d_ A(f) \right)}}  + 
\underset{(B)}{\underbrace{{\bf{d^R}} \left( \sum_{\un{n}>m} 
(-1)^{n+1} E_{m,n} (f) \right) }},\\
\end{eqnarray*}
where for the first term we have
\begin{eqnarray*}
(A)&=& (-1)^m (-1)^m d_A ^2 (f) + \sum_{\un{n}>m} (-1)^m (-1)^{n+1} E_{m,n} 
(d_A (f))\\
&=& \sum_{\un{n}>m} (-1)^{n+1} (E_{m,n} \circ (-1)^m d_A) (f),
\end{eqnarray*}and for the second term we have

\begin{eqnarray*}
(B)&=&\sum_{\un{n} >m} (-1)^{n+1} \left\{ (-1)^n d_A \left( E_{m,n} (f) 
\right) + \sum_{\un{k}>n} (-1)^{k+1} E_{n,k} \left( E_{m,n} (f) \right) 
\right\}\\
&=& \sum_{\un{n}>m} (-1)^{n+1} \left( (-1)^n d_A \circ E_{m,n}\right) (f)\\  
& &+ \sum_{\un{k} > \un{n} > m} (-1)^{n+1} (-1)^{k+1} \left( E_{n,k} \circ 
E_{m,n}\right) (f)\\
&=& \sum_{\un{n} >m} (-1)^{n+1} \left( (-1)^n d_A \circ E_{m,n} \right) (f) \\
&& + \sum_{\un{n} > \un{l}> m} (-1)^{n+1} (-1)^{l+1} \left(E_{l,n}\circ 
E_{m,l} \right) (f).
\end{eqnarray*} Thus, ${(\bf{d}^R)}^2 (f) = (A)+ (B)$ is equal to
\begin{eqnarray*}
{(\bf{d}^R)}^2 (f) &=& \sum_{\un{n}>m} (-1)^{n+1} \left\{  (E_{m,n} \circ 
(-1)^m d_A) (f) + \left( (-1)^n d_A \circ E_{m,n} \right) (f) \right. \\
&& \left. + \sum_{n > \un{l} > m} (-1)^{l+1} \left(E_{l,n}\circ E_{m,l} 
\right) (f) \right\}=0,
\end{eqnarray*}
where the last equality follows from \eqref{eqn:RCC1}.
\end{proof}

It is easy to see in the definition of a right $C$-complex 
that for $n= m+1$, $E_{m,m+1}:
A^{\bullet}_m \to A^{\bullet}_{m+1}$ is a map of chain complexes. Moreover,
the composite $E_{m+1, m+2} \circ E_{m, m+1}$ is zero in the homotopy 
category $K(\Z)$ via the homotopy $E_{m, m+2}$. Thus, a right $C$-complex 
is also a chain complex of objects of the homotopy category $K(\Z)$ of chain 
complexes.  One gets a convergent spectral sequence similar to the 
one in ~\eqref{eqn:SS1}:
\begin{equation}\label{eqn:RSS}
E^{p,q}_1  = H^q(A^{\bullet}_p) \Rightarrow H^{p+q}\left(Tot(A), 
{\bf{d}^R}\right).
\end{equation} 

Our interest in $C$-complexes is explained by the following results.
\begin{lem}\label{lem:TC-com}
Let $\sC$ be a partial dg-category. Let $A' \in {\sC}^{\oplus}$ and
let $B = \{ (B^i)_{i \in \mathbb{Z}}, q_{i,j} : B^i \to B^{j} \}$ be a 
twisted complex over $\sC$ as in Definition~\ref{defn:tComplex}. 
Assume that we have chosen distinguished subcomplexes
${\Hom_{{\sC}^{\oplus}}(A', B_m)}'$ such that the complex
$\Hom_{{\rm PreTr}(\sC)}(A', B)$ is defined as in \eqref{eqn:PreTr*}.
Then $\left(A^{\bullet}_m, d_{A_m}\right)$ is a left $C$-complex, where
$A_m = {\Hom_{{\sC}^{\oplus}}(A', B_m)}'$ and $d_{A_m}$ is its differential. 
Moreover, 
\[
\left(Tot\left(A\right), {\bf{d}^L}\right) = 
\left(\Hom_{{\rm PreTr}(\sC)}(A', B), D\right).
\]
\end{lem}
\begin{proof} We can assume that all $A', B^i \in \sC$.
%We shall denote a distinguished subcomplex
%${\Hom_{{\sC}^{\oplus}}(A', B_m)}'$ simply by 
%$\Hom_{{\sC}^{\oplus}}(A', B_m)$ assuming all the terms are defined.
Let $A_m = \Hom_{\sC}(A', B_m)'$ and $F_{m,n} = (-1)^{m+n} q_{m,n} \circ (-)$ 
for $m < n$.
Since $q_{m,n} \in \Hom^{m-n+1}_{{\sC}^{\oplus}}(B_m, B_n)$, we have
for any $f \in A^{\bullet}_m$, 
\[
F_{m,n}(f) = (-1)^{m+n} q_{m,n} \circ f \in A^{\bullet + m-n+1}_n.
\]
Moreover, the Leibniz rule 
\[
d\left(q_{m,n} \circ f\right) = d(q_{m,n}) \circ f + (-1)^{m-n+1} q_{m,n} 
\circ d(f) 
\]
for the composition in $\sC$ and ~\eqref{eqn:tComplex1} together imply
that
\[
(-1)^n d \left (q_{m,n} \circ f\right) + (-1)^{m} q_{m,n} \circ d(f)  
+ \stackrel{}{\underset {m<\un{k}<n}{\sum}} \left(q_{k,n} \circ q_{m,k} 
\circ f \right) = 0.
\]
This exactly translates to the condition ~\eqref{eqn:C-complex1} in the
definition of a left $C$-complex. This proves the first part.

For the second part, one sees from ~\eqref{eqn:PreTr*} and
~\eqref{eqn:C-complex10} that the terms of the two complexes 
$\Hom_{{\rm PreTr}(\sC)}(A, B)$ and $Tot\left(A\right)$ 
agree in each degree. Furthermore, using ~\eqref{eqn:D} and
~\eqref{eqn:C-complex2} and noting that $A'$ is a single term twisted
complex, we see that the two differentials also agree.  
\end{proof}
\begin{lem}\label{lem:TC-Rcom}
Let $\sC$ be a partial dg-category. Let $B' \in {\sC}^{\oplus}$ and
let $A = \{ (A^i)_{i \in \mathbb{Z}}, p_{i,j} : A^i \to A^{j} \}$ be a 
twisted complex over $\sC$. 
Assume that we have chosen distinguished subcomplexes 
${\Hom_{{\sC}^{\oplus}}(A^m, B')}'$ such that the complex
$\Hom_{{\rm PreTr}(\sC)}(A, B')$ is defined.
Let $
\left(B^{\bullet}_m, d_{B_m}\right) : =
\left({\Hom_{{\sC}^{\oplus}}(A^{-m}, B')}', (-1)^m d_{-m}\right),$ where $d_{-m}$ is the differential of the complex $\Hom_{\mathcal{C}^{\oplus}} (A^{-m}, B')'$. Then $\left(B^{\bullet}_m, d_{B_m}\right)$ is a right $C$-complex.
Moreover, 
\[
\left(Tot\left(B\right), {\bf{d}^R}\right) = 
\left(\Hom_{{\rm PreTr}(\sC)}(A, B'), D\right).
\]
\end{lem}
\begin{proof}
Since right $C$-complexes have not appeared before, we give a detailed proof
in this case. We first show that $\left(B^{\bullet}_m, d_{B_m}\right)$ is a
right $C$-complex. For $m < n$, Let $E_{m,n}(f) = (-1)^{\deg (f)} f \circ
p_{-n,-m}$, where $\deg (f):=r$ if  $ f \in (\Hom_{\mathcal{C}^{\oplus}} (A^{-m}, B')')^r$. Then we have
\begin{equation}\label{eqn:TC-Rcom1}
\begin{array}{lll}
\stackrel{}{\underset {m<\un{l}<n}{\sum}} (-1)^{l+1}E_{l,n}\circ E_{m,l}(f) &
= & \stackrel{}{\underset {m<\un{l}<n}{\sum}} (-1)^{l+1 + \deg (f)} E_{l,n}\left(
f \circ p_{-l,-m}\right) \\
& = & \stackrel{}{\underset {m<\un{l}<n}{\sum}} (-1)^{2\deg (f) -l+m+1}
f \circ p_{-l,-m} \circ p_{-n,-l} \\
& = & \stackrel{}{\underset {m<\un{l}<n}{\sum}} (-1)^{l+m+1}
f \circ p_{l,n} \circ p_{m,l}.
\end{array}
\end{equation}
Using the Leibniz rule, for the differential $d_{-m}$ of the complex $\Hom_{\mathcal{C}^{\oplus}} (A^{-m}, B')'$,
\[
d_{-n}\left(f \circ p_{-n,-m}\right) = (d_{-m} f) \circ p_{-n,-m} + (-1)^{\deg(f)} 
d_{-n}(p_{-n,-m}).
\] Thus, the equation ~\eqref{eqn:tComplex1} implies that
\[
(-1)^{\deg(f)}d_{-n}\left(f \circ p_{-n,-m}\right) + (-1)^{\deg(f)+1}
(d_{-m}f) \circ p_{-n,-m} 
\]
\[
\hspace*{4cm}
 + (-1)^m \stackrel{}{\underset {-n<\un{l}<-m}{\sum}} f \circ p_{l,-m}\circ p_{-n,l}
= 0. 
\]
\[
\Rightarrow 
(-1)^{\deg(f)+n}d_{B_n}\left(f \circ p_{-n,-m}\right) + (-1)^{\deg(f)+m+1} 
d_{B_m}(f) \circ p_{-n,-m} 
\]
\[
\hspace*{4cm}
+ (-1)^m \stackrel{}{\underset {-n<\un{l}<-m}{\sum}} f \circ p_{l,-m}\circ p_{-n,l}
= 0. 
\]
\[
\Rightarrow 
(-1)^{n}d_{B_n} \circ E_{m,n}(f) + (-1)^m E_{m,n} \circ d_{B_m}(f) +
\]
\[
\hspace*{4cm}
\stackrel{}{\underset {m<\un{l}<n}{\sum}}
(-1)^{l+1}E_{l,n} \circ E_{m,l} (f) = 0,
\]
where the last implication follows from the definition of $E_{m,n}$'s and
~\eqref{eqn:TC-Rcom1}.
This shows that $\left(B^{\bullet}_m, d_{B_m}\right)$ is a right $C$-complex.
The proof of the second assertion follows directly by comparing the
terms of both complexes and computing the two differentials using
~\eqref{eqn:CR-complex2} and ~\eqref{eqn:D}. Indeed, 
\begin{eqnarray*} \Hom_{\rm PreTr(\mathcal{C})} (A, B')^p &=& \bigoplus_{\un{l} + 0 + \un{m} = p}( \Hom_{\mathcal{C}^{\oplus}} (A^{-m}, B')')^l= \bigoplus_{m} (\Hom_{\mathcal{C}^{\oplus}} (A^{-m} , B')')^{p-m}\\
&=& \bigoplus_m B_m ^{p-m} = Tot(B)^p.
\end{eqnarray*}For differentials, when $f \in B_{m} ^{p-m}$, we have
\begin{eqnarray*}
\td^{\mathbf{R}} (f) &=& (-1)^m d_{B_m} (f) + \sum_{\un{n}>m} (-1)^{n+1} E_{m,n} (f) \\
&=& d_{-m} (f) + \sum_{\un{-n} < -m} (-1)^{n+1} (-1)^{p-m} f \circ p_{-n,-m},
\end{eqnarray*}while
\begin{eqnarray*}D(f) &= &(-1)^0 d_{-m} (f) + \sum_{-n} (-1)^{(p-m) + 0 + (-n) + 1} f \circ p_{-n, -m}\\
&=& d_{-m} (f) + \sum_{-n} (-1)^{p-m+n+1} f \circ p_{-n, -m} = \td ^{\mathbf{R}} (f),
\end{eqnarray*}as desired. This proves the lemma.\end{proof}

\begin{prop}\label{prop:RLtc}
Let $\mathcal{C}$ be a partial dg-category. Let 
$A = \{ (A^i)_{i\in \mathbb{Z}}, p_{i,j} : A^i \to A^j \}$, 
$B= \{ (B^i)_{i \in \mathbb{Z}}, q_{i,j} : B^i \to B^j \}$ be two twisted 
complexes over $\mathcal{C}$. Use the convention $A_i := A^{-i}$. Assume that we have chosen distinguished subcomplexes $\Hom_{\mathcal{C}^{\oplus}} (A_i, B^m )'$ for which $\Hom_{\rm PreTr(\mathcal{C})} (A, B)$ is defined. Then we 
have the following:
\begin{enumerate}
\item For fixed $i \in \mathbb{Z}$, consider the complex 
$A_{i,m} ^{\bullet} = \Hom_{\mathcal{C}^{\oplus}} (A_i, B^m)'$ and let 
$d_{A_{i,m}}$ be its differential. Then, the system $$L(A_i, B) := 
\{ (A_{i,m} ^{\bullet}= \Hom_{\mathcal{C}^{\oplus}} 
(A_i, B^m)')_{m \in \mathbb{Z}}, F_{m,n}  \}$$ is a left $C$-complex for 
some suitable $F_{m,n}$ induced from $q_{i,j}$.

\item For each $i \in \mathbb{Z}$, let $T_i = Tot ^{\cdot}(A_i, B ^{\cdot})$ 
be the total complex of the left $C$-complex $L(A_i, B)$, where the 
differential is denoted by $\td^{\mathbf{L}} _{A_i, B}$. Then, the system 
$$RL(A,B):= \{ (T_m = Tot^{\cdot} (A_m, B^{\cdot}), (-1)^m 
\td_{A_m, B} ^{\mathbf{L}})_{m \in \mathbb{Z}}, E_{m,n}\}$$ 
is a right $C$-complex for some suitable $E_{m,n}$ induced from $p_{i,j}$.

\item Let $\mathbb{T} = Tot ^{\cdot} (T_{\cdot}) = 
Tot ^{\cdot }( Tot ^{*} (A_{\cdot}, B^{*}))$ be the total complex of the 
right $C$-complex $RL (A,B)$, where the differential is denoted by 
$\td_{A,B} ^{\mathbf{RL}}.$ Then, we have
$$ (\mathbb{T}= Tot^{\cdot} (Tot^* (A_{\cdot}, B^{*})), 
\td_{A,B} ^{\mathbf{RL}}) = (\Hom_{\rm PreTr(\mathcal{C})} (A, B), D).$$ 
\end{enumerate}
\end{prop}

\begin{proof}(1) is nothing but Lemma~\ref{lem:TC-com} with $A' = A^{-i}$, 
where $F_{m,n}$ is given for $f \in A_{i,m} ^{\bullet}$ by 
$F_{m,n} = (-1)^{m+n} q_{m,n} \circ f \in A_{i,n} ^{\bullet + m - n+1}.$

(2) Let $\td= (-1)^m \td_{A_m,B} ^{\mathbf{L}}$ for simplicity. First of all, 
by the definition of the total complex $T_m$, its degree $p$-term is
\[T_m ^p = Tot^{\cdot} (A_m, B^{\cdot})^p= \bigoplus_{m ' \in \mathbb{Z}} 
A_{m,m'} ^{p-m'}= \bigoplus_{m'} \Hom_{\mathcal{C}^{\oplus}} ^{p-m'} (A^{-m}, 
B^{m'})'.\]
For $f \in A_{m, m'} ^{p-m'} \subset T_m ^p$, let 
$$E_{m,n} (f) := (-1)^{p-m'} f \circ p_{-n,-m} \in A_{n,m'} ^{p-m' -n +m +1} 
\subset T_n ^{p+m-n+1}. $$ We prove that $RL(A,B)$ is a right $C$-complex 
with respect to these $E_{m,n}$. But the perceptive reader will notice that 
when $m'\in \mathbb{Z}$ is fixed, the maps $E_{m,n}$ are defined in exactly 
same way as in Lemma~\ref{lem:TC-Rcom}, thus the relation
$$E_{m,n} \circ (-1)^m \td + (-1)^n \td \circ E_{m,n}+ 
\stackrel{}{\underset {m<\un{l}<n}{\sum}} (-1)^{l+1} E_{l,n} \circ E_{m,l} = 
0$$ 
works for all $f \in A_{m, m'} ^{p-m'}$ by the same proof. This proves (2).

(3) We prove that both $\mathbb{T}$ and 
$\Hom_{\rm PreTr(\mathcal{C})} (A, B)$ have exactly the same direct summands, and 
that on each component, $\td_{A,B} ^{\mathbf{RL}} = D$. Indeed, the degree 
$p$-term is
\begin{eqnarray*}\mathbb{T}^p &= &
\bigoplus _{m} T_{m} ^{p-m} = \bigoplus_{m} 
\bigoplus_{m'} A_{m, m'} ^{p-m-m'}\\
& =& \bigoplus_{m} \bigoplus_{m'} \Hom_{\mathcal{C}^{\oplus}} ^{p-m-m'} 
(A^{-m}, B^m)'=\Hom_{\rm PreTr(\mathcal{C})} ^p (A, B).
\end{eqnarray*} 
Regarding the differentials, let $f \in A_{m,m'} ^{p-m-m'}$, 
$\td= (-1)^m \td _{A_m,B} ^{\mathbf{L}}$ of (2), and let 
$d= d_{A_{m,m'}}$ be the differential for the complex 
$A_{m, m'} ^{\bullet} = \Hom_{\mathcal{C}^{\oplus}} (A^{-m}, B^{m'})'$. 
Then, we have 
\begin{eqnarray*}
\td_{A,B} ^{\mathbf{RL}} (f)&=& (-1)^{m} \td (f) + 
\stackrel{}{\underset {\un{n}> m}{\sum}}(-1)^{n+1} E_{m,n} (f) \\
&=& \td_{A_m, B} ^{\mathbf{L}} (f)  + \stackrel{}{\underset {\un{n}>m}{\sum}} 
(-1)^{n+1} (-1)^{p-m} f \circ p_{-n , -m}\\
&=& \left( (-1)^{m'}df + \stackrel{}{\underset {\un{n}> m}{\sum}}
F_{m', n} (f) \right) \\
&&+ \stackrel{}{\underset {\un{n}> m}{\sum}} (-1)^{n+1} (-1)^{p-m} 
f \circ p_{-n , -m}\\
&=& (-1)^{m'} df + \stackrel{}{\underset {\un{n}> m'}{\sum}}
(-1)^{m'+ n} q_{m', n}\circ f \\
&& + \stackrel{}{\underset {\un{n}> m}{\sum}} (-1)^{n+1} (-1)^{p-m} 
f \circ p_{-n , -m}.
\end{eqnarray*}
Since $A_{m,m'} ^{p-m-m'} = 
\Hom_{\mathcal{C}^{\oplus}} ^{p-m-m'} (A^{-m}, B^{m'})'$, after a suitable 
re-indexing, one immediately sees that
$\td_{A,B} ^{\mathbf{RL}} (f) = D(f).$ This finishes the proof.
\end{proof}

\begin{remk}
The Proposition~\ref{prop:RLtc} can also be stated by $(1)$ first taking the
right $C$-complexes fixing $B^j$ for each $j$, and then $(2)$ taking its associated total complexes, which along with varying $j$ form a left $C$-complex. The total complex of this out put gives the same result without affecting the final result $(3)$. We leave the detailed formulation and its proof as an exercise. 
\end{remk}

\begin{defn}\label{defn:TRC}
Let $\sC$ be a partial dg-category. We define ${\rm Tr}(\sC)$ to be 
a category such that \\
$\bullet$ $Ob\left({\rm Tr}(\sC)\right) = 
Ob\left({\rm PreTr}(\sC)\right)$ \\
$\bullet$ For any two twisted 
complexes $A = \{ (A^i)_{i\in \mathbb{Z}}, p_{i,j} : A^i \to A^j \}$
and  $B= \{ (B^i)_{i \in \mathbb{Z}}, q_{i,j} : B^i \to B^j \}$, one has 
\[\Hom_{{\rm Tr}(\sC)}(A, B) : =
H^0\left(\left(\Hom_{{\rm PreTr}(\sC)}(A,B), D\right)\right).
\]
\end{defn}
To justify the above definition, we remark that for any pair of objects
$(A,B)$ in ${\rm PreTr}(\sC)$ as above, the complex 
$\Hom_{{\rm PreTr}(\sC)}(A,B)$
depends on the choice of distinguished subcomplexes 
${\Hom_{\sC}(A_{\alpha}^i,B_{\beta}^{i'})}' \subset
\Hom_{\sC}(A_{\alpha}^i,B_{\beta}^{i'})$. If we make another choice of
the distinguished subcomplexes ${\Hom_{\sC}(A_{\alpha}^i,B_{\beta}^{i'})}''$, then by $({\rm P}5)$, we have distinguished subcomplexes ${\Hom_{\sC}(A_{\alpha}^i,B_{\beta}^{i'})}'''$ contained in both  ${\Hom_{\sC}(A_{\alpha}^i,B_{\beta}^{i'})}'$ and $ 
{\Hom_{\sC}(A_{\alpha}^i,B_{\beta}^{i'})}''$ from which we get inclusion maps
\[
{\Hom_{\sC}(A_{\alpha}^i,B_{\beta}^{i'})}' \hookleftarrow 
{\Hom_{\sC}(A_{\alpha}^i,B_{\beta}^{i'})}''' \hookrightarrow {\Hom_{\sC}(A_{\alpha}^i,B_{\beta}^{i'})}''
\]
that are quasi-isomorphisms. Then we see from Proposition~\ref{prop:RLtc}(1) that for a fixed $i \in \Z$, there is a filtered system of quasi-isomorphic
left $C$-complexes  ${\Hom_{{\sC}^{\oplus}}(A^i,B)}'$. The spectral sequence
~\eqref{eqn:SS1} and Proposition~\ref{prop:RLtc}(2) then imply that by
varying $i \in \Z$, we get a filtered system of quasi-isomorphic
right $C$-complexes  ${\Hom_{{\rm PreTr}(\sC)}(A^i,B)}'$. Finally,
Proposition~\ref{prop:RLtc}(3) and spectral sequence ~\eqref{eqn:RSS} together
imply that there is a well-defined filtered system of quasi-isomorphic
complexes ${\Hom_{{\rm PreTr}(\sC)}(A,B)}'$. In particular,
$H^0\left(\Hom_{{\rm PreTr}(\sC)}(A,B), D\right)$ is canonically defined.
In particular, ${\rm Tr}(\sC)$ is an honest category.

If $\sD$ is a partial full subcategory of ${\rm PreTr}(\sC)$ in the sense of 
Definition~\ref{defn:PreTRC}, then we define ${\rm Tr}(\sD)$ to be the
category whose objects are same as those of $\sD$ and whose morphisms
are defined as in ${\rm Tr}(\sC)$. It easily follows from the above discussion
that ${\rm Tr}(\sD)$ is an honest category and in fact is a genuine
full subcategory of ${\rm Tr}(\sC)$.  

We shall call ${\rm Tr}(\sC)$ to be the \emph{homotopy} category of the partial
category $\sC$. This terminology is inspired by the example of the
dg-category $C(R)$ of complexes of modules over a commutative ring $R$, where
the category ${\rm Tr}(C(R))$ is indeed the usual homotopy category $K(R)$
of complexes of $R$-modules.
\begin{prop}\label{prop:TriangC}
Let $\sC$ be a partial dg-category. Then the homotopy category ${\rm Tr}(\sC)$
is a triangulated category. Moreover, any functor $t : \sC \to \sD$ of
partial dg-categories induces an exact functor ${\rm Tr}(t) :
{\rm Tr}(\sC) \to {\rm Tr}(\sD)$ of triangulated categories.
\end{prop} 
\begin{proof}
We only describe the shift functor and the distinguished triangles in
${\rm Tr}(\sC)$. The rest of the proof follows exactly like 
\cite[Section~4]{Ha1}. We skip the details and refer to \emph{ibid.} to see that
all axioms of a triangulated category are satisfied.

Let $A = \{ (A^i)_{i\in \mathbb{Z}}, q_{i,j} : A^i \to A^j \}$
be an object of ${\rm Tr}(\sC)$. The shift functor $A \mapsto A[1]$  
is given by ${A[1]}^i = A^{i+1}$ and $q[1]_{i,j} = (-1)^{i+j+1}
q_{i+1,j+1}$. For a morphism $u : A \to B$, $u[1]$ is given by
${(u[1])}^{i,j} = (-1)^{i+j} u^{i+1,j+1}$.

The cone of a morphism $u = (u^{i,j}) : A = 
(A^i, q_{i,j}) \to (B^{i'}, r_{i',j'}) = B$ is defined as an object $C =
(C^k, t_{k,l})$, where
\[
C^i = A^{i+1} \oplus B^i
\]
and

\[
t_{i,j} = \left(  \begin{array}{clcr}
(-1)^{i+j+1}q_{i+1,j+1} & 0 \\
u^{i+1,j} & (-1)^{i+j}r_{i,j}
\end{array}
\right).
\]
There are natural morphisms $\alpha (u) : B \to C$ given by
\[
\alpha (u)^{i,j} =
\left(  \begin{array}{clcr}
0 \\
(-1)^i \delta_{i,j} 1_{L^i}
\end{array} \right)
: B^i \to A^{i+1} \oplus L^i
\]
and $\beta (u)^{i,j} : C \to A[1]$ is given by
\[
\beta (u)^{i,j} = \left(\delta_{i,j} 1_{A^{i+1}} \right) :
A^{i+1} \oplus B^i \to A^{i+1},
\]
where $\delta_{i,j}$ is $0$ if $i\not = j$, $1$ if $i=j$, so that there is a distinguished triangle
\[
A \to B \to C \to A[1]
\]
in ${\rm Tr}(\sC)$. Such triangles are called the standard distinguished
triangles and an arbitrary distinguished triangle in ${\rm Tr}(\sC)$ is the
one isomorphic to a standard one.
\end{proof}

\begin{remk} If $\sD$ is a partial full subcategory of
${\rm PreTr}(\sC)$ in the sense of Definition~\ref{defn:PreTRC}, then we have 
seen before that ${\rm Tr}(\sD)$ is a full subcategory of ${\rm Tr}(\sD)$.
However, this may not be the inclusion of triangulated categories since
${\rm Tr}(\sD)$ may not be closed under the cone construction. 
For example, one could take $\sD$ as those twisted complexes in ${\rm PreTr}
(\sC)$ in Section~\ref{subsection:PdgC} where $q_{i,j}$'s are only higher
Chow cycles. It is then easy to see that ${\rm Tr}(\sD)$ is not a 
triangulated subcategory of ${\rm Tr}(\sC)$.  
\end{remk} 

%\begin{exercise}[{\emph{cf.} \cite[Proposition 2]{BK}}] Show that ${\rm PreTr} (-)$ is a closed operation in the following sense: for any partial dg-category $\mathcal{C}$, we have a natural functor
%$${\rm PreTr} ({\rm PreTr}(\mathcal{C})) \to {\rm PreTr} (\mathcal{C})$$ whose $H^0 : {\rm Tr} ({\rm PreTr}(\mathcal{C})) \to {\rm Tr} (\mathcal{C})$ is an equivalence of triangulated categories.
%\end{exercise}

\section{Additive cycle complexes and some properties}
\label{section:cyclecomplex}
In this section, we review the definition of additive cycle complexes
from \cite{KP} and also study some of their properties which we shall need
in this paper. We begin with a recall of the cubical version of
Bloch's higher Chow complexes from  \cite[p. 298]{Levine}. 
Set $\P^1:=\Proj\, k[Y_0,Y_1]$, and set $\square^n:=(\P^1\setminus\{1\})^n$. We use the coordinates $(y_1, \cdots, y_n )$ for $\square^n$. A \emph{face} $F \subset \square^n$ is a closed subscheme defined by equations of the form $\{ y_{i_1} = \epsilon_1, \cdots, y_{i_s} = \epsilon_s\}$, where each $\epsilon_j$ is $0$ or $\infty$. For each $\epsilon = 0, \infty$ and each $i= 1, \cdots, n$, we let $\iota_{n, i, \epsilon} : \square^{n-1} \to \square^n$ be the closed immersion given by $(y_1, \cdots, y_{n-1}) \mapsto (y_1 , \cdots, y_{i-1} , \epsilon, y_i, \cdots, y_{n-1})$. The schemes $\iota_{n, i, \epsilon} (\square^{n-1})$ are called the codimension $1$ faces of $\square^n$.

\begin{defn}Let $X$ be a $k$-variety. Let $\underline{z} ^q (X, n)$ be the free abelian group generated by closed irreducible subvarieties $Z \subset X \times \square^n$ that intersect all faces of $\square^n$ properly, i.e. in the right codimensions. The cycle $\iota_{n, i, \epsilon}^* (Z) \in \underline{z} ^q (X,n-1)$ is denoted by $\partial_i ^\epsilon (Z)$. Define the boundary map as $\partial := \sum_{i=1} ^n (-1)^i (\partial_i ^{\infty} - \partial_i ^0)$.

Let $\underline{z} ^q (X, n)_{\dgn}$ be the subgroup of $\underline{z} ^q (X, n)$ given by the degenerate cycles, i.e. cycles obtained by pulling back via the projections
$X \times \square^n \to X \times \square^{n-1}.$ Define $z^q (X, n):= \underline{z} ^q (X, n)/\underline{z} ^q (X, n)_{\dgn}.$ One checks that the boundary map $\partial$ descends to $z^q (X, n)$, and $\partial^2 = 0$. This is the cubical higher Chow complex of $X$, and its homology is the higher Chow group denoted by $\CH^q (X, n)$.\end{defn}

\subsection{Additive cycle complexes}\label{subsection:cyclecomplex1}
We follow the notations of \cite{KP} to define the additive cycle complexes.
For a $k$-scheme $V$,
let $V^N$ be the normalization of $V_{\rm red}$. Set $\A^1:=\Spec k[t]$, $\G_m:=\Spec k[t,t^{-1}]$. For $n \ge 1$, let $B_n = {\G}_m \times \square^{n-1}$, 
$\ov{B}_n =  {\A}^{1} \times ({\P}^{1})^{n-1}$
and $\widehat{B}_n = {\P}^{1} \times ({\P}^{1})^{n-1}$, with the coordinates $(t, y_1, \cdots , y_{n-1})$ on $\widehat{B}_n$.

Let $F^1_{n,i}$, for $i=1,\ldots, n-1$, be the divisor $\{ y_i = 1 \}$, and $F_{n,0}$
the divisor $\{t=0\}$ on $\widehat{B}_n$. Let $F^1_n:=\sum_{i=1} ^{n-1} F^1_{n,i}$ on
$\widehat{B}_n$. A {\em face} $F$ of $B_n$ is defined by equations of the form $y_{i_1}=\epsilon_1, \ldots,  y_{i_s}=\epsilon_s$ where each $\epsilon_j$ is $0$ or $\infty$. For each $\epsilon=0, \infty$ and each $i = 1, \cdots, n-1$, let $\iota_{n,i,\epsilon}: B_{n-1}\to B_n$ be the inclusion $(t,y_1,\ldots, y_{n-2}) \mapsto (t,y_1,\ldots, y_{i-1}, 
\epsilon,y_i,\ldots, y_{n-2})$, that gives a codimension $1$ face.

\subsubsection{Modulus conditions}\label{subsection:MD}
The additive higher Chow cycles satisfy one additional property, other than the proper intersection with faces, called the \emph{modulus condition}. We consider two such conditions for which the moving lemma of \cite{KP} works, which is essential in this paper:
\begin{defn}\label{defn:M5}
Let $X$ be a $k$-variety, and let $V$ be an integral closed subvariety 
of $X \times B_n$. Let $\ov V$ be the Zariski closure of $V$ in $X \times
\widehat{B}_n$, and let $\nu : {\ov V}^N \to  X \times {\widehat {B}}_n$ be 
the normalization of $\ov V$. Fix an integer $m \ge 1$. 
\begin{enumerate}
\item
We say that $V$ satisfies the modulus $m$ condition $M_{sum}$ on $X \times B_n$, if as Weil divisors on 
${\ov V}^N$, we have $$
(m+1)[{\nu}^*(F_{n,0})] \le [\nu^*(F^1_n)].$$
\item  
We say that $V$ satisfies the modulus $m$ condition $M_{ssup}$ on $X \times B_n$, if there exists an integer $1 \le i \le n-1$ such that $$
(m+1)[{\nu}^*(F_{n,0})] \le [{\nu}^*(F^1_{n,i})]$$
as Weil divisors on ${\ov V}^N$.
\end{enumerate}
\end{defn}
We often  say that $V$ has the modulus condition $M$ without mentioning $m$.

\begin{defn}\label{defn:AdditiveComplex}
Let $M$ be the modulus condition $M_{sum}$ or $M_{ssup}$.
Let $X$ be an equi-dimensional $k$-variety, and let $r,m$ be integers with 
$m\ge1$. \\
\\
(0) $\un{\TZ}_r(X, 1; m)_M$  is the free abelian group 
on integral closed subschemes $Z$ of $X \times {\G}_m$ of dimension $r$.\\
\\
For $n>1$,  $\un{\TZ}_r(X, n; m)_M$ is the free abelian group 
on integral closed subschemes $Z$ of $X \times B_n$ of dimension $r+n-1$ 
such that: \\ \\
\noindent (1) For each face $F$ of $B_n$, $Z$ intersects  $X \times F$
properly on $X \times B_n$.

\noindent (2) $Z$ satisfies the modulus $m$ condition $M$ on $X \times B_n$.   
\end{defn}
If $d= \dim X$, we write for $q \ge 0$
\[
\un{\TZ}^q(X, n; m)_M = \un{\TZ}_{d+1-q}(X, n; m)_M.
\]

As shown in \cite{KP}, one can check that if $Z \subset X \times B_n$ 
satisfies the above conditions $(1)$ and $(2)$, then every component of ${\iota_{n,i,\epsilon}}^*(Z)$ also satisfies these conditions on $X \times B_{n-1}$. As before, we let $\un{\TZ}^q(X, n; m)_{M, \dgn}$ be the subgroup generated by the degenerate cycles.
\begin{defn}
The {\em additive higher Chow complex} $\TZ^q(X, \bullet ; m)_M$ of $X$ in codimension $q$ and with modulus $m$ condition $M$ is the non-degenerate complex
$$
\TZ^q(X, \bullet; m)_M: = \un{\TZ}^q(X, \bullet; m)_M
/{\un{\TZ}^q(X, \bullet; m)_{M, \dgn}}.$$
The boundary map is
$\partial = \sum_{i=1} ^{n-1} (-1)^i (\partial^{\infty}_i
- \partial^0_i)$. It satisfies $\partial ^2 = 0$. 
The homology
$
\TH^q(X, n; m)_M: = H_n (\TZ^q(X, \bullet ; m)_M)$ for $n \ge 1
$
is the 
{\it additive higher Chow group}  of $X$ with modulus $m$ condition $M$. 
\end{defn}
We shall drop the subscript $M$ from the notations. All results of this paper work for both the modulus conditions.

\subsubsection{Total higher Chow complex}

This paper deals with the higher Chow cycles and the additive higher Chow cycles altogether:

\begin{defn}
The \emph{total higher Chow complex} of $X$ of codimension $q$ with respect to modulus $m \geq 1$ is the direct sum of complexes
$$z^q (X, \bullet; m) := z^q (X, \bullet) \oplus  \TZ^q (X, \bullet; m).$$ Its degree $n$ homology will be denoted by 
$$\CH^q (X, n; m) := \CH^q (X, n) \oplus  \TH^q (X, n;m).$$
By convention, for $m=0$, we let $z^q (X, \bullet; 0)$ be the higher Chow 
complex $z^q (X, \bullet)$, and let $\CH^q (X, n; 0)$ be the higher Chow 
group $\CH^q (X, n)$.
\end{defn}

We shall need the following functoriality properties of the cycle complexes.
\begin{lem}[Push-forward and pull-back]\label{lem:functoriality-lemma}Let $X, Y, Z, X', Y'$ be $k$-varieties.
\begin{enumerate} 
\item If $f: X \to Y$ is a projective morphism, then the push-forward
$f_* : z^q (X, \bullet ; m) \to z^{q'}( Y, \bullet ; m),$ $q':=q + { \dim} Y - { \dim} X $, is well-defined on the level of complexes. 

If $f: X \to Y$ and $g: Y \to Z$ are two projective morphisms, then we have $(g\circ f)_* = g _* \circ f_*$.

\item If $f: X \to Y$ is a flat morphism, then the pull-back
$f^* : z^q (Y, \bullet; m) \to z^q (X, \bullet ; m)$ is well-defined on the level of complexes. 

If $f: X \to Y$ and $g: Y \to Z$ are two flat morphisms, then we have $(g \circ f)^* = f^* \circ g^*$.
\item If we have a Cartesian square 
$$ \begin{CD} X' @>{g'}>> X \\
@V{f'}VV @VV{f}V\\
Y' @>>{g}> Y\end{CD}$$ where $f$ is flat and $g$ is projective, then, as maps 
of complexes $z^q (Y', \bullet ; m) \to z^{q'}(X, \bullet; m)$, we have $f^* \circ g_* = g'_* \circ {f'} ^*,$ where $q'= q + \dim X - \dim X' = q + \dim Y - \dim Y'$.

\end{enumerate}
\end{lem} 
\begin{proof}(1), (2) for the higher Chow part follow from the cubical version of \cite[Proposition 1.3]{Bl1}. For the additive part, \cite[Lemmas 3.6, 3.7]{KL} prove them for a slightly different modulus condition, so we state here that the same works for the modulus conditions $M_{sum}$ and $M_{ssup}$ as follows: (1) is a special case of a stronger statement \cite[Proposition 5.2]{KP} that if $f: X \to Y$ is projective, and $Z \subset X \times B_n$ is an admissible additive cycle, then its projective image $(f \times 1_{B_n}) (Z)$ is an admissible additive cycle. The push-forward is the projective image when the map $Z \to (f \times 1_{B_n}) (Z)$ is generically finite, and $0$ if not. For (2), the proof in \cite[Lemma 3.7]{KL} works without change.

(3) This statement is part of more general compatibility of projective 
push-forward and flat pull-back of any cycles in Cartesian squares as in 
\cite{Fulton}. \end{proof}

\subsection{Operations of correspondences on total cycles}
Total higher Chow cycles on $X \times Y$ can induce two important partially defined operations on cycles in total higher Chow complexes. 
These operations are induced by the following external products of cycles.

\subsubsection{External product} Given two $k$-varieties $X, Y$, we have the external products
$$\boxtimes : z^{q_1} (X, n_1;m) \otimes z^{q_2} (Y, n_2; m) \to z^{q} (X \times Y, n; m)$$
where $q = q_1 + q_2, n = n_1 + n_2$, for each integral closed admissible subschemes $Z_1 \subset X \times \square^{n_1}$ or $X \times B_{n_1}$ and $Z_2 \subset Y \times \square^{n_2}$ or $Y \times B_{n_2}$, given by $Z_1 \boxtimes Z_2 = \tau _* (Z_1 \times Z_2)$, where $$\tau: \tuborg X \times \square^{n_1} \times Y  \times \square^{n_2} \to X \times Y \times \square^{n_1} \times \square^{n_2} \\ X \times \square^{n_1} \times  Y \times \mathbb{G}_m \times \square^{n_2 -1} \to X \times Y \times \mathbb{G}_m \times \square^{n_1} \times \square^{n_2 -1}\\
 X \times \mathbb{G}_m \times \square^{n_1 -1} \times Y \times \square^{n_2} \to X \times Y \times \mathbb{G}_m \times \square^{n_1 -1} \times \square^{n_2} \sluttuborg$$
 are the corresponding transpositions. By convention, the product $\boxtimes$ of two additive admissible cycles is zero. If only one of $Z_1$ and $Z_2$ is an additive admissible cycle, then $Z_1 \boxtimes Z_2$ is an additive admissible cycle by \cite[Lemma~4.2]{KL}. Hence the above external products $\boxtimes$ is well-defined. Note that $\boxtimes$ is  distributive over sums.

\subsubsection{Cup product} Given $X \in \SmProj/k$, we have partially defined products
\begin{eqnarray*}\cup_X  : z^{q_1} (X, n_1; m) \otimes z^{q_2} (X, n_2; m) \dto z^{q} (X, n; m)\end{eqnarray*}
where $q= q_1 + q_2, n = n_1 + n_2$, given by the formula $$ Z_1 \cup_X Z_2 = \delta_X ^* (Z_1 \boxtimes Z_2), $$ if the pull-back via the diagonal $\delta_X : X \to X \times X$ makes sense. While $\boxtimes$ is always defined, the pull-back $\delta_X ^*$ is defined only on a distinguished subcomplex $z^{q} (X \times X, \bullet; m)'$ by 
Lemma~\ref{lem:functoriality-lemma2}-(1).
\begin{remk}The notation $\cap_X$ was used in \cite[\S 4]{KL} for the second and the third cases of the above. We uniformly use the notation $\cup_X$ for simplicity.
\end{remk}

\begin{lem}\label{lem:projection-formula}Let $X, Y \in \SmProj/k$, and let $f: X \to Y$ be a morphism of $k$-varieties.
\begin{enumerate}
\item [(0)] The cup product is associative, and distributive over sums, whenever they are defined.
\item We have
$f^* ( Z_1 \cup_Y Z_2) = f^* (Z_1) \cup_X f^* (Z_2),$ whenever all expressions are defined.
\item If $f$ is both flat and projective, then whenever all expressions are define, we have the projection formulas
$$ f_* (f^* (Z_1) \cup_X Z_2) = Z_1 \cup_Y f_* (Z_2),  f_* (Z_1 \cup_X f^* (Z_2)) = f_* (Z_1) \cup_Y Z_2.$$
\end{enumerate}
\end{lem}
\begin{proof}(0) Consider for each $i=1,2,3$, integral closed admissible subschemes $Z_i \subset X \times \times \square^{n_i}$ or $X \times B_{n_i}$.
We show that $(Z_1 \cup_X Z_2) \cup_X Z_3 = Z_1 \cup_X (Z_2 \cup_X Z_3)$, if all cup-products are defined. Consider the following commutative diagram
$$\begin{CD}
X \times X @>{\delta_X \times 1}>> X \times X \times X \\
@A{\delta_X}AA @AA{1 \times \delta_X} A\\
X @>>{ \delta_X}> X \times X\end{CD}$$ from which we get
$\delta_X ^* ( \delta_X \times 1)^* = \delta_X ^* (1 \times \delta_X^*)$ by 
Lemma~\ref{lem:functoriality-lemma}-(2). Since $\boxtimes$ is associative, we have
\begin{eqnarray*}
(Z_1 \cup _X Z_2) \cup Z_3 &=& \delta_X ^*( \delta_X ^* ( Z_1 \boxtimes Z_2) \boxtimes Z_3)\\
&=& \delta_X ^* ((\delta_X \times 1)^* ((Z_1 \boxtimes Z_2) \boxtimes Z_3) \\
&=& \delta_X ^* ((1 \times \delta_X)^* (Z_1 \boxtimes (Z_2 \boxtimes Z_3))\\
&=& Z_1 \cup_X (Z_2 \cup_X Z_3).
\end{eqnarray*}This proves the associativity. The distributive law is obvious by definition.

(1) Note that we have a commutative diagram
$$\begin{CD}
X @>{\delta_X}>> X \times X \\
@V{f}VV @VV{f \times f}V \\
Y @>>{\delta_Y}> Y \times Y
\end{CD}$$from which we get 
\begin{eqnarray*}
f^* \circ \delta_Y ^* &= &(\delta_Y \circ f)^*  \ \ \ \ \ \ \ \  \mbox{ 
(by Lemma~\ref{lem:functoriality-lemma}-(2))} \\
&=& ((f \times f)\circ \delta_X)^* \\
&=& \delta_X ^* \circ (f \times f)^*.
\end{eqnarray*}
Hence, by a direct calculation we have
\begin{eqnarray*}
f^* (Z_1 \cup_Y Z_2)& =& f^* (\delta_Y ^* (Z_1 \boxtimes Z_2)) \\
&=& \delta_X ^* ( (f\times f)^* (Z_1 \boxtimes Z_2))\\
&=& \delta_X ^* (f ^* (Z_1) \boxtimes f^* (Z_2))
\end{eqnarray*} as desired. For (2), one can follow 
\cite[Theorem~4.10]{KL}. \end{proof}

\subsubsection{Push-forward by correspondences}
\label{subsubsection:pushforward} Let $X, Y \in \SmProj/k$, and let $v \in z^{q_2} (X \times Y, n_2;m)$. Then, we have partially defined push-forward maps of complexes
\begin{eqnarray*}v_* :  z^{q_1} (X, \bullet;m) \dto  z^{q } (Y, \bullet + n_2;m),
\end{eqnarray*}where $q = q_1 + q_2 - \dim X$, given by $v_* (Z): = p_{Y*}  \left( v \cup_{X \times Y} {p_X ^{*}} (Z) \right),$ and $p_X, p_Y $ are the obvious projections. If one writes $v = (\alpha, f)$ where $\alpha$ is the higher Chow cycle and $f$ is the additive cycle, then one has $v_* = \alpha_* + f_*$.

\subsubsection{Composition by correspondences}\label{Hanamura partially defined composition} Let $X, Y, Z  \in \SmProj/k$. We have partially defined compositions

\begin{eqnarray*} (-) \circ (-) : z^{q_2} (Y \times Z, \bullet; m) \otimes z^{q_1} (X \times Y, \bullet ; m) \dto z^{q} (X \times Z, \bullet ; m),
\end{eqnarray*}
where $q= q_1 + q_2 - \dim Y$, given by
$$ v \otimes u \mapsto v \circ u := {p_{XZ*} ^{XYZ}} \left( {p_{YZ}^{XYZ*}}  (v) \cup {p_{XY}^{XYZ*}}  (u) \right),$$ where $\cup = \cup_{X \times Y \times Z}$, and $p_{XZ} ^{XYZ}$, etc. are the obvious projections. Since $\cup$ is distributive over sums of cycles, if one writes a cycle $v$ as $v= (\alpha,f)$, where $\alpha_i$ is a higher Chow cycle and $g_i$, an additive one, we deduce the composition law
\begin{equation}\label{eqn:total-composition-in-moving}(\alpha_2, f_2) \circ (\alpha_1, f_1) = (\alpha_2 \circ \alpha_1, \alpha_2 \circ f_1 + f_2 \circ \alpha_1).\end{equation}

\begin{lem}\label{lem:associativity-for-composition}Let $X, Y, Z , W\in \SmProj/k$.
\begin{enumerate}
\item For three higher Chow cycles $\alpha_1$ on $X \times Y$, $\alpha_2$ on $Y \times Z$, and $\alpha_3$ on $Z \times W$, we have $(\alpha_3 \circ \alpha_2) \circ \alpha_1 = \alpha_3 \circ (\alpha_2 \circ \alpha_1),$ if all compositions are defined.

\item Let $f \in \TZ ^{q_1} (X \times Y, n_1;m)$, $\alpha_1 \in z^{q_2} (Y \times Z, n_2)$, $\alpha_2 \in z^{q_3} (Z \times W, n_3)$. Then we have 
$ (\alpha_2 \circ \alpha_1 ) \circ f = \alpha_2 \circ (\alpha_1 \circ f),$ if all compositions are defined.

Similarly, for cycles on appropriate spaces, we have
$ (\alpha_1 \circ f) \circ \alpha_2= \alpha_1 \circ (f \circ \alpha_2),$ and
$ (f \circ \alpha_1) \circ \alpha_2 = f \circ (\alpha_1 \circ \alpha_2)$ if all compositions are defined, where $\alpha_i$ are higher Chow cycles, and $g$ are additive higher Chow cycles.
\item The composition law of ~\eqref{eqn:total-composition-in-moving} for the total correspondences is associative whenever the compositions are defined.
\end{enumerate}
\end{lem}

\begin{proof}(1) is proven in \cite{Ha1}. 

(2) We let $p^{XYZ} _{XY}$, etc. be the obvious projections, but the projections from $ X \times Y \times Z \times W$ to, say $X \times Z$, will be denoted by $p_{XZ}$ instead of $p_{XZ} ^{XYZW}$.  We prove the first equation $ (\alpha_2 \circ \alpha_1 ) \circ f = \alpha_2 \circ (\alpha_1 \circ f).$ From the RHS, 
we have
\begin{eqnarray*}
&&\alpha_2 \circ (\alpha_1 \circ f) \\
&=&{ p_{XW*}^{XZW}} \left\{ {p_{ZW} ^{XZW*}}(\alpha_2) \cup {p_{XZ} ^{XZW*} }( \alpha_1 \circ f)\right\}\\
&=& { p_{XW*}^{XZW}} \left[ {p_{ZW} ^{XZW*}} (\alpha_2) \cup p_{XZ} ^{XZW*} \left\{ p_{XZ*} ^{XYZ}\left( p_{YZ} ^{XYZ*} (\alpha_1) \cup p_{XY} ^{XYZ^*} (f)   \right) \right\}  \right] \\
&=& { p_{XW*}^{XZW}} \left[ {p_{ZW} ^{XZW*}} (\alpha_2) \cup p_{XZW*} ^{} \left\{ p_{XYZ} ^{*} \left(p_{YZ} ^{XYZ*} (\alpha_1) \cup p_{XY} ^{XYZ^*} (f)   \right) \right\}  \right] \\
&& ( \mbox{by Lemma~\ref{lem:functoriality-lemma}-(3)})\\
&=& { p_{XW*}^{XZW}} \left[ {p_{ZW} ^{XZW*}} (\alpha_2) \cup p_{XZW*} ^{} \left\{p_{YZ} ^{*} (\alpha_1) \cup p_{XY} ^{*} (f) \right\} \right]\\
&& (\mbox{by Lemma~\ref{lem:projection-formula}-(1) and 
Lemma~\ref{lem:functoriality-lemma}-(2)})\\
&=& { p_{XW*}^{XZW}} \left[ p_{XZW*} ^{} \left\{  p _{XZW} ^{*} \left({p_{ZW} ^{XZW*}} (\alpha_2) \right) \cup \left( p_{YZ} ^{*} (\alpha_1) \cup p_{XY} ^{*} (f) \right) \right\} \right]\\
&& (\mbox{by the projection formula, Lemma~\ref{lem:projection-formula}-(2)})\\
&=&p_{XW*} \left\{ p_{ZW} ^* (\alpha_2) \cup \left( p_{YZ} ^* (\alpha_1) \cup p_{XY} ^* (f) \right) \right\} \\
&=& p_{XW*} \left\{ \left( p_{ZW} ^* (\alpha_2) \cup p_{YZ} ^* (\alpha_1) \right) \cup p_{XY} ^* (f) \right\} (\mbox{by Lemma~\ref{lem:projection-formula}(0)).}
\end{eqnarray*}
By the same kind of calculations with LHS, we get to the last expression. Hence, we get first equation. The other two equations are proven in exactly the same fashion.

(3) Let  $v_i =(\alpha_i, f_i)$, $i=1, 2, 3$ be total correspondences for which $v_1 \circ v_2$, $v_2 \circ v_3$,  and $(v_1 \circ v_2) \circ v_3$, $v_1 \circ (v_2 \circ v_3)$ are defined. Note that

\begin{eqnarray*}
(v_1 \circ v_2) \circ v_3&=& ((\alpha_1, f_1) \circ (\alpha_2, f_2)) \circ (\alpha_3, f_3) \\ & =&  (\alpha_1 \circ \alpha_2, \alpha_1 \circ f_2 + f_1 \circ \alpha_2) \circ (\alpha_3, f_3) \\
&= &((\alpha_1 \circ \alpha_2) \circ \alpha_3, (\alpha_1 \circ \alpha_2 )\circ f_3\\
&& + (\alpha_1 \circ f_2) \circ \alpha_3 + (f_1 \circ \alpha_2) \circ \alpha_3),\\
v_1 \circ (v_2 \circ v_3) &=& (\alpha_1 , f_1) \circ ((\alpha_2 , f_2) \circ (\alpha_3, f_3)) \\
&=& (\alpha_1, f_1) \circ (\alpha_2 \circ \alpha_3, \alpha_2 \circ f_3 + f_2 \circ \alpha_3) \\
&=& (\alpha_1 \circ (\alpha_2 \circ \alpha_3), \alpha_1  \circ (\alpha_2 \circ f_3)\\
&& + \alpha_1 \circ( f_2 \circ \alpha_3) + f_1 \circ (\alpha_2 \circ \alpha_3)). 
\end{eqnarray*} Thus, (1) and (2) imply the equality $(v_1 \circ v_2) \circ v_3 = v_1 \circ (v_2 \circ v_3).$
\end{proof}

\section{Moving lemma and distinguished subcomplexes}\label{section:ML}
In this section, we define a class of distinguished subcomplexes for the 
total higher Chow complexes, and study its properties. This is technically 
the most important part. The new ingredient behind this definition is the 
moving lemma for additive higher Chow groups of smooth projective varieties 
in \cite{KP} and its refinement discussed below. In this section, $m \ge 0$ is a fixed 
integer and all the results hold for any of the
modulus conditions considered in Section~\ref{section:cyclecomplex}. 

\subsection{A refined moving lemma}\label{subsection:Rmoving}
The moving lemma for higher Chow groups from \cite{Bl1} 
(\emph{cf.} \cite{Ha1} for the cubical version) and additive higher Chow 
groups from \cite[Theorem 4.1]{KP} of smooth projective varieties together 
imply the following form of moving lemma:
\begin{thm}\label{thm:KPmoving}
Let $X \in \SmProj/k$. Let 
$\mathcal{W}$ be a finite set of irreducible locally closed subsets of $X$. 
Then, the inclusion
$z_{\mathcal{W}} ^q (X, \bullet; m) \hookrightarrow z ^q (X, \bullet ; m)$ is 
a quasi-isomorphism.
\end{thm}
\begin{remk}\label{remk:moving-e} By \cite[Remarks 4.3, 4.4]{KP}, the above 
theorem is equivalent to that the inclusion 
$z^q _{\mathcal{W}, e} (X, \bullet;m) \hookrightarrow z^q (X, \bullet ; m)$
is a quasi-isomorphism for all set functions 
$e : \mathcal{W} \to \mathbb{Z}_{\geq 0},$ where 
$$z ^q _{\mathcal{W},e} (X, \bullet ; m) := 
z^q _{\mathcal{W}, e} (X, \bullet) \oplus  
\TZ^q _{\mathcal{W}, e} (X, \bullet ; m),$$ and 
${z} ^q _{\mathcal{W}, e} (X, n )$ 
(resp. $\TZ^q _{\mathcal{W} , e} (X, n;m)$) is defined as follows: first, let 
$\underline{z}^{q}_{\sW,e }(X, n)$ 
(resp. $\underline{\TZ} ^q _{\mathcal{W}, e}(X, n;m)$) be the subgroup of 
$\underline{z}^{q}(X, n)$ (resp. $\underline{\TZ} ^q (X, n;m)$) generated by 
integral closed subschemes
$Z \subset X \times \square^n$ (resp. $Z \subset X \times B_n$) such that 
$ \codim _{W \times F} (Z \cap (W \times F)) \geq q - e(\mathcal{W})$ for all 
$W \in \mathcal{W}$ and all faces $F$ of 
$\square^n$ (resp. all faces $F$ of $B_n$).
We let ${z}^{q}_{\sW,e}(X, \bullet ;m) $ be the image of 
$\underline{z}^{q}_{\sW,e}(X, \bullet ;m)= \un{z} ^q _{\mathcal{W}, e}(X, n) \oplus \un{\TZ} ^q _{\mathcal{W},e} (X, \bullet ; )$ in ${z}^{q}(X, \bullet ;m) $ via the 
projection modulo the degenerate cycles.
\end{remk}

This paper requires a bit finer form of moving lemma than 
Theorem~\ref{thm:KPmoving}. We allow the following more general collections 
$\mathcal{W}$ of varieties:

\begin{defn}[c.f. {\cite[Definition 2.1]{KL}} 
{\cite[p.112]{Ha1}}]\label{defn:extended-case} Let $X \in \SmProj/k$ and let 
$T_1, \cdots, T_n$ be finitely many $k$-schemes of finite type over $k$. Let 
$\mathcal{W}$ be a finite set of irreducible locally closed subsets 
$W_i \subset X \times T_i$ for $ i = 1, \cdots, N.$ For each face 
$F  \subset \square^n$ and $F \subset B_n$, let $p_{F, i} : 
X \times F \times T_i \to X \times T_i$ be the projection. 

Let $\underline{z} ^q _{\mathcal{W}} (X, \bullet ; m) \subset 
\underline{z} ^q (X, \bullet ; m)$ be the direct sum of subcomplexes 
$\underline{z} ^q _{\mathcal{W}} (X, \bullet)$ in $\un{z} ^q (X, \bullet)$ 
and $\underline{\TZ}^q _{\mathcal{W}} (X, \bullet; m)$ in 
$\un{\TZ} ^q (X, \bullet; m)$,  where 
$\underline{z} ^q _{\mathcal{W}} (X, \bullet)$ is generated by integral 
closed subschemes $Z \subset X \times \square^n$ such that, additionally, for each face 
$F \subset \square^n$, two sets $p_{F, i} ^{-1} (W_i)$ and 
$(Z \cap (X \times F)) \times T_i$ intersect properly on 
$X \times F \times T_i$ for all $i = 1, \cdots, N$. The complex 
$\underline{\TZ} ^q _{\mathcal{W}} (X, \bullet; m)$ is defined similarly.

The image of $\underline{z} ^q  _{\mathcal{W}}(X, \bullet; m)$  in 
$z ^q (X, \bullet ; m)$, under the projection modulo degenerate cycles,  is 
called a \emph{distinguished subcomplex} of $z^q (X, \bullet ; m)$. 
Similarly, one defines the complexes $z^q _{\mathcal{W}} (X, \bullet)$ and 
$\TZ^q _{\mathcal{W}} (X, \bullet; m)$ modulo degenerate cycles, and they are 
called distinguished subcomplexes of $z^q (X, \bullet)$ and 
$\TZ^q (X, \bullet;m)$, respectively. If the reference to the set 
$\mathcal{W}$ is not necessary, then we simply write 
$z^q (X, \bullet;m)'$ for any distinguished subcomplex.
\end{defn}

\begin{remk}\label{remk:traditional-case}As a special case, take all 
$T_i = \Spec (k)$ for $i = 1, \cdots, N$. Then, we recover the complex 
$z^q _{\mathcal{W}} (X, \bullet; m)$ in Theorem~\ref{thm:KPmoving}.
\end{remk}

We now discuss a refined version of the moving lemma:

\begin{thm}\label{thm:refined-moving} Let $X\in \SmProj/k$, 
and let $\mathcal{W}$ be a finite set of irreducible $k$-varieties as in 
Definition~\ref{defn:extended-case}. Then, the inclusion
$z_{\mathcal{W}} ^q (X, \bullet; m) \hookrightarrow z ^q (X, \bullet ; m)$ is 
a quasi-isomorphism.
\end{thm}
\begin{proof}One can prove it using arguments similar to those in 
\cite[Proposition~2.2]{KL} together with Theorem~\ref{thm:KPmoving}, and 
Remark~\ref{remk:moving-e}: for each $W_i \subset X \times T_i$, $i = 1, 
\cdots, N$, in the set $\mathcal{W}$, we form the constructible subsets of $X$
$$C_{i,d} := \{ x \in X | (x \times T_i) \cap W_i 
\mbox{ contains a component of dimension } \geq d \}.$$ 
Write each $C_{i,d} \backslash C_{i, d-1}$ as a union of irreducible 
locally closed subsets $C_{i,d} ^j$. Let 
$\mathcal{C}:= \{ C_{i, d} ^j | i, d, j \}$, and let 
$e: \mathcal{C} \to \mathbb{Z}_{\geq 0}$ be the set-theoretic function 
defined by $e(C_{i,d} ^j ) := \dim W_i - d - \dim C_{i,d} ^j,$ which is 
always $\geq 0$. One can then check by comparing the defining conditions,
that 
$$\underline{z} ^q _{\mathcal{C}, e} (X, \bullet ; m) = 
\underline{z} ^q _{\mathcal{W}} (X, \bullet ; m).$$ 
By the moving lemma of Theorem~\ref{thm:KPmoving} and 
Remark~\ref{remk:moving-e}, the image of the left is quasi-isomorphic to 
$z^q (X, \bullet; m)$, thus so does the image of the right hand side. This 
proves the theorem. 
\end{proof}
The following obvious result is very frequently used in this paper, so we 
record it here. (\emph{cf.} $({\rm P}5)$ in Definition \ref{defn:Pdg}, Proposition \ref{prop:PartialCC})

\begin{lem}\label{lem:finite-intersection}For $X \in \SmProj/k$, let $z^q (X, \bullet; m)', z^q (X, \bullet; m)''$ be two distinguished subcomplexes. Then, there exists a distinguished subcomplex $z^q (X, \bullet; m)'''$ contained in both $z^q (X, \bullet; m)'$ and $z^q (X, \bullet; m) ''$.
\end{lem}

\begin{proof}Let $\mathcal{W}'$, $\mathcal{W}''$ be the finite sets as in 
Definition~\ref{defn:extended-case} that give the complexes $z^q (X, \bullet; m)', z^q (X, \bullet; m)''$, respectively.

Take all the set $T_i$'s used to give $\mathcal{W}'$ and $\mathcal{W}''$ (as in Definition~\ref{defn:extended-case}), and collect all of $W_i$'s in $\mathcal{W}'$ and $\mathcal{W}''$ to define $\mathcal{W}$. This gives
$z^q (X, \bullet; m)''':= z^q _{\mathcal{W}} (X, \bullet; m) = z^q (X, \bullet; m)' \cap z^q( X, \bullet; m)''.$
\end{proof}
The following useful lemma is backed by the refined moving lemma:
\begin{lem}\label{lem:functoriality-lemma2}Let $X, Y,Z $ be $k$-varieties.
\begin{enumerate}
\item Suppose $Y\in \SmProj/k$, and let $f: X \to Y$ be any morphism. Then, there exits a distinguished subcomplex $z^q (Y, \bullet ; m)'$ on which the pull-back $f^* : z^q (Y, \bullet ;m)' \to z^q (X, \bullet ; m)$ is well-defined on the level of complexes. 

\item Let $X, Y \in \SmProj/k$, and let $f: X \to Y$ be any morphism. Then, given any distinguished subcomplex $z^q (X, \bullet;m)'$, there exists a distinguished subcomplex $z^q (Y, \bullet; m)'$ on which the pull-back $f^*$ is well-defined, and we have
$f^* \left( z^q (Y, \bullet ; m)' \right) \subset z^q (X, \bullet; m)'.$

\item Let $X, Y \in \SmProj/k$, and let $f: X \to Y$ be a projective morphism. Then, given any distinguished subcomplex $z^{q'} (Y, \bullet; m)'$, where $q' := q + \dim Y - \dim X$, there exists a distinguished subcomplex $z^q (X, \bullet; m)'$ such that
$f_* \left( z^q (X, \bullet; m)' \right) \subset z^{q'} (Y, \bullet; m)'.$
\end{enumerate}
\end{lem}

\begin{proof} A similar but weaker statement was proven in 
\cite[Theorem~7.1]{KP} using Theorem~\ref{thm:KPmoving}. A more efficient 
proof is provided here with Theorem~\ref{thm:refined-moving}. (1) As in 
Definition~\ref{defn:extended-case}, we take $T = X$, and take $W = {}^t \Gamma_f \subset Y \times T$, the transpose of the graph of $f$. Take $z^q (Y, \bullet ; m) ':= z^q _{\{ {}^t \Gamma_f \}} (Y, \bullet ; m)$. Then, it gives a natural pull-back
$p_Y ^* : z^q _{ \{ {}^t \Gamma_f \}} (Y, \bullet; m) 
\to z^q _{\{ {}^t \Gamma_f\}} (Y \times X, \bullet ; m),$ where $p_Y$ is the
projection $X \times Y \to Y$, and the subscript $\{ {}^t \Gamma_f\}$ is in the sense of Remark~\ref{remk:traditional-case}. Composing with the Gysin chain map induced by the regular embedding ${}^t {\rm Graph}_f : X \to Y \times X$ (see \cite[Corollary~7.2]{KP}) $ z^q _{\{ {}^t \Gamma_f \}} (Y \times X, \bullet, m) \to z^q (X, \bullet; m),$ one gets $f^* : z^q _{\{ {}^t \Gamma_f \}} (Y, \bullet ; m) \to z^q (X, \bullet ; m),$ as desired.

(2) Let $\mathcal{W}$ be a set of $W_i \subset X \times T_i$, for some $k$-schemes with $i = 1, \cdots, N$ with the desired properties as in 
Definition~\ref{defn:extended-case} that gives the given distinguished subcomplex $z^q(X, \bullet; m)'$.

Then one takes for $\mathcal{W}'$, the collection of the sets ${p_Y} ^{-1} (W_i) \subset Y \times S_i$, with $S_i:= X \times T_i$ for $i=1, \cdots, N$, and ${}^t \Gamma_f\subset Y \times S_{N+1}$ with $S_{N+1} := X$, where $p_Y: Y \times X \to Y$ is the projection. Take $z^q (Y, \bullet ; m)':= z^q _{\mathcal{W}'} (Y, \bullet; m)$. Then
$f^* \left( z^q _{\mathcal{W}'} (Y, \bullet; m) \right) 
\subset z^q _{\mathcal{W}} (X, \bullet ; m) = z^q (X, \bullet; m)'$ as desired. 

(3) We drop the codimensions when no confusion arises to simplify our notations. Let $\mathcal{W}$ be such that $z _{\mathcal{W}} (Y, \bullet ; m)$ is the given distinguished subcomplex $z_{\mathcal{W}} (Y, \bullet; m)'$. Assume that the set $\mathcal{W}$ is given by the irreducible closed subvarieties $W_i \subset Y \times T_i$, $i = 1, \cdots, N$, for some $k$-schemes $T_i$. 

Consider $(f \times 1_{T_i})^{-1} (W_i) \subset X \times T_i$, and write $W_{ij}$ for the irreducible components of $(f \times 1_{T_i})^{-1} (W_i)$. Let $\mathcal{W}' = \{ W_{ij}| i,j\}$. Then we have $f_* \left( z_{\mathcal{W}'} (X, \bullet; m)\right) \subset z_{\mathcal{W}} (Y, \bullet; m)=z (Y, \bullet ; m)'.$ This finishes the proof.
\end{proof}

\subsubsection{Distinguished subcomplexes and the operations}Let's have a closer look at the above partially defined operations using the refined moving lemma. The following proposition summarizes some essential results we need later in the paper. This generalizes \cite[Propositions~1.4, 1.5]{Ha1} 
(\emph{cf.} \cite[Proposition~2.5]{KL}) to include additive higher Chow cycles.

\begin{prop}\label{prop:grand-summary}Let $X, Y, Z, W \in \SmProj/k$. Let $q_i, n_i \geq 0$ be integers. Then we have the following properties:

\begin{enumerate}
\item [$(1a)$] Given $w \in z^{q_2} (X, n_2;m)$, there exists a distinguished subcomplex $z^{q_1} (X, \bullet; m)'$ on which $w \cup_X (-)$ is defined, and similarly for $(-) \cup_X w$.

\item [$(1b)$] Given $w \in z^{q_2} (X, n_2;m)$ and a given distinguished subcomplex $z^{q} (X, \bullet; m)'$, with $q= q_1 + q_2$, there exists a distinguished subcomplex $z^{q_1} (X, \bullet; m)'$ on which $w \cup _X (-)$ is defined, and we have
$w \cup_X \left(z^{q_1} (X, \bullet; m)'\right) \subset z^{q} (X, \bullet+ n_2; m)'.$ Similarly for $(-) \cup_X w$.

\item [$(2a)$] Given $v \in z^{q_2} (X \times Y, n_2)$, there exists a distinguished subcomplex $z^{q_1} (X , \bullet;m) '$ on which $v_*$ is defined.

\item [$(2b)$] In addition to $(2a)$, given $w \in z^{q_3} (Y \times Z, n_3)$ such that $\beta \circ v$ is defined, there exist distinguished subcomplexes $z^{q_1} (X, \bullet; m)'$  and $z^{q} (Y, \bullet; m) '$, where $q = q_1 + q_2 - \dim X$, such that
\begin{enumerate}
\item [${\rm (i)}$] both $v_*$ and $(w\circ v)_*$ are defined on $z^{q_1} (X, \bullet; m)',$
\item [${\rm (ii)}$] $w_*$ is defined on $z^{q} (Y, \bullet; m)'$,
\item [${\rm (iii)}$] $v_* \left( z^{q_1} (X, \bullet; m)' \right) \subset z^{q} (Y, \bullet+ n_2 
; m)',$ and
\item [${\rm (iv)}$] $w_* \circ v_* = (w \circ v)_*$ on $z^{q_1} (X, \bullet ;m)'$.

\end{enumerate}

\item [$(3a)$] Given $v \in z^{q_2} (Y \times Z, n_2 ;m)$, there exists a distinguished subcomplex $z^{q_1} (X \times Y, \bullet; m)'$ on which 
$ v \circ (-)$ is defined.

\item [$(3b)$] In addition to $(3a)$, given $w \in z^{q_3} (Z \times W, n_3;m)$ such that $w \circ v$ is defined, there exists distinguished subcomplexes $ z^{q_1} (X \times Y, \bullet;m)'$ and $z^{q} (X \times Z, \bullet; m)'$, where $q = q_1 + q_2 - \dim Y$, such that
\begin{enumerate}
\item [${\rm (i)}$] $v \circ (-), (w \circ v) \circ (-)$ are defined on $z^{q_1} (X \times Y, \bullet; m)'$, 
\item [${\rm (ii)}$] $w \circ (-)$ is defined on $z^{q} (X \times Z, \bullet; m)'$, 
\item [${\rm (iii)}$] $v \circ \left( z^{q_1} (X  \times Y, \bullet; m)' \right) \subset z^{q} (X \times Z, \bullet + n_2; m) '$, and
\item [${\rm (iv)}$] $w \circ \left( v \circ (-) \right) = ( w \circ v) \circ (-)$ on $z^{q_1} (X \times Y, \bullet;m)'$.
\end{enumerate}
The same works for compositions from the right.
\item [$(3c)$] The same works for any finite sequence of the above operations.
\end{enumerate}
\end{prop}
\begin{proof}
$(1a)$ is a special case of $(1b)$. For $(1b)$, given $w\in z^{q_2} (X, n_2;m)$ and a distinguished subcomplex $z^{q} (X, \bullet ; m)'$, by 
Lemma~\ref{lem:functoriality-lemma2}-(2) there exists a distinguished subcomplex $z^{q} _{\mathcal{W}_1} (X \times X, \bullet ; m)$ on which $\delta_X ^*$ is defined, and
$ \delta_X ^* (z^{q} _{\mathcal{W}_1} (X \times X, \bullet + n_2 ; m))$ is contained in $z^{q}  (X, \bullet + n_2; m)'.$ Then, it is enough to find a set $\mathcal{W}_2$ for which we have
$w \boxtimes z^{q_1} _{\mathcal{W}_2} (X, \bullet; m) \subset z^{q} _{\mathcal{W}_1} (X \times X , \bullet + n_2; m).$ We may assume $w$ is an irreducible closed subvariety of $X \times \square^{n_2}$ by 
Lemma~\ref{lem:finite-intersection}. If $\mathcal{W}_1$ is given by 
$\mathcal{W}_1 = \{ W_i \subset X \times X \times T_i | i = 1, \cdots, N \},$ then for $\mathcal{W}_2$ we take
$\mathcal{W}_2 = \{ W_i \subset X \times S_i | i = 1, \cdots, N \} \cup \{ w \subset X \times S_{N+1} \},$ where $S_i = X \times T_i$ for $i =1, \cdots, N$ and $S_{N+1} = \square^{n_2}$. This proves $(1b)$.

$(2a)$ Recall that for $Z \in z^{q_1} (X, \bullet;m)$, if defined, the push-forward $v_* (Z)$ is given by the expression $v_* (Z): = {p_{Y*} ^{XY}} \left( v \cup_{X \times Y} {p_X ^{XY*}} (Z) \right).$ 

By $(1a)$, there exists a distinguished subcomplex $z^{q_1} (X \times Y, \bullet; m)'$ on which the product $v \cup_{X \times Y} (-)$ is defined on $z^{q_1} (X \times Y, \bullet; m)'$. By Lemma~\ref{lem:functoriality-lemma2}-(2), there exists a distinguished subcomplex $z^{q_1} (X, \bullet; m)'$ such that $p_X^{XY*}  \left( z^{q_1} (X, \bullet; m)' \right) \subset z^{q_1} (X \times Y, \bullet; m)'$. Since $p_{Y} ^{XY}$ is projective, $p_{Y*} ^{XY}$ is everywhere defined on $z^{q_1} (X \times X, \bullet; m)'$. Hence everywhere on $z^{q_1} (X, \bullet; m)'$, the push-forward $v_*$ is defined, proving $(2a)$.

$(2b)$ We first show ${\rm (iv)}$ that $w_* \circ v_* = (w \circ v)_*$, if defined.  The notations $p_{X} ^{XY}$ etc. are the obvious projections, while the projections from $X\times Y \times Z$ are simply denoted by $p_{YZ}$, instead of $p_{YZ} ^{XYZ}$. For a cycle $Z$ on $X$ for which $(w \circ v)_* (Z)$, $(w_* \circ v_*) (Z)$ are defined, we have
\begin{eqnarray*}
(w \circ v) _* (Z) &=& p_{Z*} ^{XZ} \left( (w \circ v) \cup p_X ^{XZ*} (Z) \right) \\
&=& p_{Z*} ^{XZ} \left\{ p_{XZ*} ^{XYZ} \left( p_{YZ} ^{XYZ*}( w) \cup  p_{XY} ^{XYZ*} (v) \right) \cup p_X ^{XZ*} (Z) \right\} \\
&=& p_{Z*} ^{XZ} \left[ p_{XZ*}  \left\{ \left( p_{YZ} ^{*} (w) \cup p_{XY} ^{*} (v) \right) \cup p_{XZ} ^{*} \left( p_X ^{XZ*} (Z) \right) \right\} \right] \\
&&:\mbox{by projection formula, Lemma~\ref{lem:projection-formula}(2)}\\
&=& p_{Z*}  \left\{ \left( p_{YZ} ^{*} (w) \cup p_{XY} ^{*} (v) \right) \cup p_X ^{*} (Z) \right\} :\mbox{Lemma~\ref{lem:functoriality-lemma}(1,2)}\\
&=& p_{Z*}  \left\{ p_{YZ} ^{*} (w) \cup \left( p_{XY} ^{*} (v) \cup p_X ^{*} (Z) \right) \right\}:\mbox{Lemma~\ref{lem:projection-formula}(0)}\end{eqnarray*}
\begin{eqnarray*}
{ \ \ \ \ \ \ \ \ \ }&=& p_{Z*} ^{YZ} \left[ p_{YZ*}  \left\{ p_{YZ} ^{*} (w) \cup \left( p_{XY} ^{*} (v) \cup p_X ^{*} (Z) \right) \right\} \right]:\mbox{Lemma~\ref{lem:functoriality-lemma}(1)}\\
&=& p_{Z*} ^{YZ} \left\{ w \cup p_{YZ*}  \left( p_{XY} ^{*} (v) \cup p_X ^{*} (Z) \right) \right\}:\mbox{Lemma~\ref{lem:projection-formula}(2)}\\
&=& p_{Z*} ^{YZ} \left\{ w \cup p_{YZ*}  \left( p_{XY} ^{*} (v) \cup p_{XY} ^{*} p_X ^{XY*} (Z) \right) \right\}:
\mbox{Lemma~\ref{lem:functoriality-lemma}(2)}\\
&=& p_{Z*} ^{YZ} \left[ w \cup p_{YZ*} \left\{ p_{XY} ^{*} \left(v \cup p_X ^{XY*} (Z) \right)\right\} \right]:\mbox{Lemma~\ref{lem:functoriality-lemma}(1)}\\
&=& p_{Z*} ^{YZ} \left[ w \cup p_{Y} ^{YZ*} \left\{ p_{Y*} ^{XY}\left( v \cup p_X ^{XY*} (Z) \right) \right\} \right]:\mbox{Lemma~\ref{lem:functoriality-lemma}(3)}\\
&=& w_* \left\{ p_{Y*} ^{XY} \left( v \cup p_X^{XY*} (Z) \right) \right\} = (w_* \circ v_* )(Z).
\end{eqnarray*}
The rest of $(2b)$ is similar to $(2a)$, but a bit more is involved. Whenever the codimensions we consider are apparent, we will drop them for simplicity. Given $w\in z (Y \times Z,n_3;m)$, by $(2a)$ we have a distinguished subcomplex $z(Y, \bullet; m)'$ on which $w_*$ is defined, which gives ${\rm (ii)}$. By 
Lemma~\ref{lem:functoriality-lemma2}-(3), one can find a distinguished subcomplex $z (X \times Y, \bullet; m)'$ such that 
\begin{equation}\label{eqn:aa}
p_{Y*} ^{XY} (z (X \times Y, \bullet+ n_2; m)') \subset z (Y, \bullet+ n_2; m)'.\end{equation}
Now, for $v \in z (X \times Y, n_2;m)$, by $(1b)$ we have a distinguished subcomplex $z (X \times Y, \bullet; m)'$ on which $v \cup_{X\times Y} (-)$ is well-defined, and 
\begin{equation}\label{eqn:bb}
v \cup_{X \times Y} \left( z (X \times Y, \bullet; m)' \right) \subset z (X \times Y, \bullet + n_2; m)'.\end{equation} 
By Lemma~\ref{lem:functoriality-lemma2}-(2), one can find a distinguished subcomplex $z (X, \bullet; m)'$ such that 
\begin{equation}\label{eqn:cc}p_X ^{XY*} (z (X, \bullet; m)') \subset z (X \times Y, \bullet; m)'.\end{equation}
On the other hand, by $(2a)$ applied to $v_*$ and $(w\circ v)_*$ (with Lemma~\ref{lem:finite-intersection}), one can replace $z (X, \bullet; m)'$ by a smaller distinguished subcomplex, denoted by the same symbols, $z (X, \bullet; m)'$ on which $v_*$ and $(w \circ v)_*$ are all defined, which gives  ${\rm (i)}$. 

Now combining the above, one sees that 
\begin{eqnarray*}v_* ( z (X, \bullet; m)')& = & p_{Y*} ^{XY} \left( v \cup p_{X} ^{XY*} ( z (X, \bullet; m)') \right)\\
&\subset & p_{Y*} ^{XY} \left( v \cup z (X \times Y, \bullet; m)'\right) \mbox{~\eqref{eqn:cc}}\\
&\subset & p_{Y*} ^{XY} \left(z(X \times Y, \bullet + n_2; m)'\right) \mbox{\eqref{eqn:bb}}\\
&\subset & z (Y, \bullet + n_2; m)' \mbox{~\eqref{eqn:aa}},
\end{eqnarray*}which proves ${\rm (iii)}$. This proves all of $(2b)$.

$(3a)$ Again, we drop the codimensions from our notations, whenever no confusion arises. We let $p_{XY}$, etc. be the projections from $X \times Y \to Z$ to $X \times Y$, etc. Given $v \in z (Y \times Z, n_2;m)$, for the fixed $p_{YZ} ^* (v) \in z(X \times Y \times Z, n_2;m)$, by (1a), there exists a distinguished subcomplex $z(X \times Y \times Z, \bullet, m)'$ on which the operation $p_{YZ} ^* (v) \cup (-)$ is well-defined. Now, by 
Lemma~\ref{lem:functoriality-lemma2}-(2), one can find a distinguished subcomplex $z(X \times Y, \bullet; m)'$ such that $p_{XY} ^* \left( z(X \times Y, \bullet; m)' \right) \subset z (X \times Y \times Z, \bullet; m)'.$ Since $p_{XZ*}$ is everywhere defined, on $ z(X \times Y, \bullet; m)'$ the operation $ v \circ (-) = p_{XZ*} \left( p_{YZ} ^* (v) \cup p_{XY} ^* (-) \right)$ is well-defined. This solves $(3a)$.

$(3b)$ The part ${\rm (vi)}$ is Lemma~\ref{lem:associativity-for-composition}-(3). The rest of the proof is similar to $(3a)$, but it is a bit more involved. In fact, one can imitate the arguments for $(2b)$. The reader is encouraged to try its proof following $(2b)$ with suitable changes. The arguments for the compositions from the right are similar. $(3c)$ is obvious from all of the above.
\end{proof}

\begin{cor}\label{cor:additive-associativity-1}For $X_i, Y_i \in \SmProj/k$, given finitely many $v_i \in z^{s_i }( Y_i \times Y_{i+1}, n_i;m)$, $i = 1, \cdots, N$, and $w_j \in z^{s'_j} (X_{j+1} \times X_j, n'_j;m)$, $j= 1, \cdots , N'$ for which the compositions $v_N \circ \cdots  \circ v_1$ and $w_1 \circ \cdots \circ w_{N'}$ are defined, there exists a distinguished subcomplex $z^{q} (X_1 \times Y_1, \bullet ; m)' $ on which the composition 
$v_N \circ \cdots \circ v_1 \circ ( -) \circ w_1 \circ \cdots \circ w_{N'} $ is well-defined without ambiguity.
\end{cor}

\begin{proof}To find a distinguished subcomplex on which the composition is well-defined, one repeatedly applies Proposition~\ref{prop:grand-summary}, and 
Lemma~\ref{lem:finite-intersection}. Once the compositions are defined, by the associativity the composition in question is unambiguous. This finishes the proof.
\end{proof}

\section{The category ${\sD}{\sM}(k;m)$}\label{section:DM}
In this section, we construct our category of mixed motives 
${\sD}{\sM}(k;m)$ over $k[t]/{(t^{m+1})}$ using the results of 
the previous sections. The strategy is to define a ``category''
$\sC$ which is shown to be a partial dg-category using the results of
Section~\ref{section:ML}. The desired category ${\sD}{\sM}(k;m)$ will be the
pseudo-abelian hull of the homotopy category ${\rm Tr}(\sC)$ of $\sC$. So we 
first describe our partial dg-category. We fix an integer
$m \ge 0$.

\subsection{Partial dg-category $\sC$}\label{subsection:PdgC}
The partial dg-category $\mathcal{C}$ has for objects,
the pairs $(X, r)$ for $X \in \SmProj/k$ and $r \in \mathbb{Z}$. The objects
have the product structure via $(X, r) \otimes (Y, s) = (X \times Y, r+s)$, 
the dual structure $(X,r)^{\vee} = (X, \dim X - r)$, and the internal hom 
$\underline{\sHom} \left((X,r) , (Y,s)\right) = (X,r)^{\vee} \otimes (Y,s)$. 
If $m \ge 1$, one associates for each object, a right bounded complex 
\begin{equation}\label{eqn:Pdcom*}
\mathcal{Z}\left((X, r); m\right) := z^r (X, - \bullet; m) = 
z^r (X, - \bullet) \oplus \TZ^r (X, - \bullet ; m),
\end{equation} 
which is the direct sum of the higher Chow complex and the additive higher 
Chow complex, seen as a cohomological complex by using $- \bullet$. 
For $m = 0$, by convention
\begin{equation}\label{eqn:Pdcom*1}
\mathcal{Z}\left((X, r); 0\right) := z^r (X, - \bullet; 0) = 
z^r (X, - \bullet),
\end{equation}
the higher Chow complex of $X$.
For two 
objects $(X, r), (Y, s) \in {\rm Ob} (\mathcal{C})$, one defines the morphism 
to be the above complex for the internal hom 
$\underline{\sHom}\left((X,r), (Y,s)\right)$, namely,
$$\hom_{\mathcal{C}} ((X,r), (Y,s)) = 
\mathcal{Z} \left((X,r)^{\vee} \otimes (Y,s); m \right).$$ 

Given two cycles $v_i = (\alpha_i, f_i)$, $i=1,2$ where 
$\alpha_1, \alpha_2$ are higher Chow cycles on $X\times Y$ and $Y\times Z$ 
respectively and $f_1,f_2$ are additive higher Chow cycles on $X \times Y$ 
and $Y \times Z$ respectively, we defined their composition by
\begin{equation}\label{eqn:composition}
v_2 \circ v_1 = (\alpha_2, f_2) \circ (\alpha_1, f_1) : = 
(\alpha_2 \circ \alpha_1, \alpha_2 \circ f_1 + f_2 \circ \alpha_1)
\end{equation}
whenever these compositions of cycles are defined.

If $(X, r)$ is an object of $\sC$, we define the unit endomorphism
as the morphism of chain complexes
\begin{equation}\label{eqn:unitM}
\mathbb{I}_{(X,r)} :\Z \to \sZ\left((X,r);m\right)
\end{equation}
given by $1 \mapsto [{\Delta}_X]$, where $[\Delta_X]$ is the class of the
diagonal in $z^{\dim X}(X \times_k X, 0)$. Note that
\begin{eqnarray*} \Hom_{\mathcal{C}} ((X,r), (X,r))& =& 
\mathcal{Z} ((X, \dim X - r) \otimes (X, r);m) \\
&=& \mathcal{Z} (X\times X, \dim X;m )\\
&=& z ^{\dim X} (X \times X, - \bullet) 
\oplus \TZ^{\dim X} (X \times X, - \bullet;m).
\end{eqnarray*}
Since $\TZ^{\dim X} (X \times X, 0;m) = 0$, we see that 
$\mathbb{I}_{(X,r)}$ is a well defined map of complexes. It is well-known
and easy to check that for any $X, Y \in \SmProj/k$,
the compositions
\begin{equation}\label{eqn:unitM*}
\Delta_Y \circ (-) : z^r (X \times Y, \bullet ; m) \dto z^r 
(X \times Y, \bullet ; m), 
\end{equation}
\[ 
(-) \circ \Delta_X : z^r (X \times Y, \bullet ; m) \dto 
z^r ( X \times Y, \bullet ; m) 
\]
are partially defined morphisms which are identity.

\begin{prop}\label{prop:PartialCC} 
$\sC$ is a partial dg-category. 
\end{prop}
\begin{proof}
To prove the proposition, we first describe the classes of \emph{distinguished
subcomplexes} and then verify all the axioms of Definition~\ref{defn:Pdg}.
For $A = (X,r), B = (Y, s) \in \sC$, we let the class $S\left(A,B\right)$ to be the
class of distinguished subcomplexes $z^{r+s} _{\mathcal{W}}\left(X \times Y, -\bullet ;
m\right)$ in the sense of Definition~\ref{defn:extended-case}.

The ``unit'' endomorphism of axiom $({\rm P}3)$, 
$\mathbb{I}_{(X,r)} :\Z \to \sZ\left((X,r);m\right)$  
is defined in ~\eqref{eqn:unitM}. The axiom  $({\rm P}4)$ is verified
in Proposition~\ref{prop:grand-summary}.  The axiom  $({\rm P}5)$ follows
directly from the stronger assertion in Lemma~\ref{lem:finite-intersection}. 
The associativity part
of axiom $({\rm P}6)$ is proven in Proposition~\ref{prop:grand-summary}
$(3b){\rm (iv)}$, and the identity action of the unit morphism is shown in
~\eqref{eqn:unitM*}. Thus $\sC$ is a partial dg-category.
\end{proof}
\subsection{The category ${\sD}{\sM}(k;m)$}\label{section:DM1}
Let $\sC$ be the partial dg-category described in Section ~\ref{subsection:PdgC}.
Proposition~\ref{prop:PartialCC} implies that $\sC$ is indeed a partial 
dg-category. In particular,
it follows from Proposition~\ref{prop:TriangC} that ${\rm Tr}(\sC)$ is a
triangulated category.
\begin{defn}\label{defn:DM*}
We define ${\sD}{\sM}(k;m)$ to be the pseudo-abelian hull of the triangulated
category ${\rm Tr}(\sC)$. 
\end{defn}
It follows from \cite[Theorem 1.5]{BS} that ${\sD}{\sM}(k;m)$ is also a
triangulated category such that there is an exact inclusion
functor $\tau : {\rm Tr}(\sC) \inj {\sD}{\sM}(k;m)$. The category 
${\sD}{\sM}(k;m)$ will be called the triangulated category of mixed motives
over the ring $k[t]/{(t^{m+1})}$ for given $m \ge 0$. 

\begin{remk}\label{remk:HDM}
It is easy to see from the definition of $\sC$ and from ~\eqref{eqn:Pdcom*1}
that ${\sD}{\sM}(k;0)$ is the same as the integral version of Hanamura's triangulated category of
mixed motives ${\sD}{\sM}(k)$ over $k$.
\end{remk} 

\begin{defn}We define the \emph{motive functor with the modulus $m$ augmentation} $h : \SmProj/k \to {\sD}{\sM}(k;m)$ as
\begin{equation}\label{eqn:mfunctor}
h(X) = \left(\left(X, 0\right), 0\right),
\end{equation}where the morphisms $f: X \to Y$ are sent to the graph $\Gamma_f \in \CH^{\dim Y} (X \times Y, 0)$.
Note that the term $((X,0),0)$ on the right of \eqref{eqn:mfunctor} is a twisted complex where the correspondences
$q_{i,j}$'s are all zero.
\end{defn}

\subsection{Some structural properties of $\mathcal{DM}(k;m)$}
\label{subsection:structure}
We now discuss some structural properties of ${\sD}{\sM}(k;m)$ which 
essentially follow from our construction and the proofs of similar results in
\cite{Ha1}.

\subsubsection{Duals} For an object $A = (A_i, q_{i,j}) \in {\rm Ob} ({\rm PreTr}(\mathcal{C}))$, define its dual object $A^{\vee}= ((A^{\vee})_i, (q^{\vee})_{i,j})$ by the relations

$$(A^{\vee})_i := (A_{-i})^{\vee}, \ \ \mbox{where the RHS ${\vee}$ is in the sense of duals for } \mathcal{C^{\oplus}},$$
$$ (q^{\vee})_{i,j} := (-1)^{ij-j+1} \  {}^t q_{-j,-i}.$$
That $A^{\vee}$ is an object of ${\rm PreTr}(\mathcal{C})$ can be checked easily.

\subsubsection{Monoidal structure}For two objects 
$A = (A_i, q_{i,j})$, $A' = (A'_i, q' _{i,j})$, we define their product 
$A \otimes A' = (M_i, h_{i,j})$ to be given as follows: for each $i$, we let

$$M_i := \bigoplus_{i_1 + i_2= i} A_{i_1} \otimes A_{i_2} ', \ \ \mbox{where the RHS $\otimes$ is for $\mathcal{C^{\oplus}}$}.$$ If $A_{i} =\bigoplus_{\alpha} A_{i, \alpha}$ $A' _i = \bigoplus_{\beta} A' _{i, \beta}$, where $A_{i, \alpha}, A' _{i, \beta} \in \mathcal{C}$, then we can write $M_i$ as
$$M_i = \bigoplus _{i_1 + i_2 = i} \bigoplus_{ \alpha, \beta} A_{i_1, \alpha} \otimes A_{i_2, \beta}.$$

The morphism $h_{i,j} : M_i \to M_j$ of degree $i-j+1$ in $\mathcal{C}^{\oplus}$ is given by combining various morphisms $$ h^{(i_1, i_2, \alpha, \beta)} _{(j_1, j_2, \alpha', \beta')}:  A_{i_1, \alpha} \otimes A_{i_2, \beta}' \to A_{j_1, \alpha'} \otimes A_{j_2, \beta'}', \ \ i_1 + i_2 = i, j_1 + j_2 = j.$$ Here, they are defined as follows:
\begin{eqnarray*} &&\mbox{$(a)$ if $\beta = \beta'$ and $i_2 = j_2$, then}\\
&& h^{(i_1, i_2, \alpha, \beta)} _{(j_1, j_2, \alpha', \beta')}:= (-1)^{i_2 (j_1-i_1-1)} q_{i_1, j_1, \alpha, \alpha'} \otimes 1,\\
&&\mbox{where $q_{i_1, j_1, \alpha, \alpha'}: A_{i_1, \alpha} \to A_{j_1, \alpha'}$ is the map given from $q_{i,j}$}, \\
&&\mbox{$(b)$ if $\alpha = \alpha'$ and $i_1 = j_1$, then}\\
&& h^{(i_1, i_2, \alpha, \beta)} _{(j_1, j_2, \alpha', \beta')}:= (-1)^{i_1} 1 \otimes q'_{i_2, j_2, \beta, \beta'},\\
&&\mbox{where $q'_{i_2, j_2, \beta, \beta'} : A'_{i_2, \beta} \to A'_{j_2, \beta'}$ is the map given from $q'_{i,j}$},\\
&& \mbox{$(c)$ for all other cases, we let $h^{(i_1, i_2, \alpha, \beta)} _{(j_1, j_2, \alpha', \beta')} =0 $.}
 \end{eqnarray*} This system gives an object of ${\rm PreTr}(\mathcal{C})$.

\subsubsection{Internal homs} For two objects $A, A' \in {\rm PreTr} (\mathcal{C})$, define the internal hom by $$\underline{\sHom} (A, A') := A ^{\vee} \otimes A'.$$

The above three operations thus give objects in ${\sD}{\sM}(k;m)$.
\begin{defn}{(\emph{Unit object})}\label{defn:unitob}
The object $\un{\Z} = \un{\Z}(0) \in {\rm Tr}(\mathcal{C})$ is defined by 
$A=(A_i, q_{i,j})$, where
$$A_i = \tuborg \left( \Spec (k), 0 \right), &\mbox{ if } i = 0 , \\ 0, & \mbox{ if } i \not = 0,\sluttuborg$$ and $q_{i,j}= 0$ for all $i, j$.
\end{defn}
\begin{prop}[{\cite[p. 140]{Ha1}}]\label{prop:tensor-properties}
For objects $A, A', A''$ of ${\sD}{\sM}(k;m)$, we have
\begin{enumerate}
\item \emph{Associativity:} $(A \otimes A') \otimes A'' = A\otimes (A' \otimes A'')$.
\item \emph{Unit object:} $\un{\mathbb{Z} } \otimes A = A \otimes \un{\mathbb{Z}} = A$.
\item \emph{Product and dual:} $(A \otimes A') ^{\vee} = {A'}^{\vee} \otimes A^{\vee}$.
\item \emph{Product and hom:} There are functorial isomorphisms  
$${\rm adj} : \Hom_{{\sD}{\sM}(k;m)}(A'', A^{\vee} \otimes A') = 
\Hom_{{\sD}{\sM}(k;m)} (A'' \otimes A, A').$$
\item \emph{Reflexivity:} There are functorial isomorphisms
$$i_A: A \to A^{\vee \vee},$$ given by $(-1)^i$ on $A_i$.
\end{enumerate}
\end{prop}
\begin{proof} See \emph{loc. cit}.
\end{proof}

\subsection{The category of mixed Tate motives}\label{subsection:Ttwist} 
For each $n \in \mathbb{Z}$, the \emph{Tate objects} $\un{\Z}(n)$ in
${\sD}{\sM}(k;m)$ are defined as
$$\un{\Z}(n) := \left( \Spec (k), n)\right) [ -2n],$$ 
i.e. $\left( \Spec (k), n\right)$ is in degree $2n$. 

The \emph{$n$-th Tate twist} $(-) (n): {\mathcal{DM}}(k;m) \to 
\mathcal{DM}(k;m)$ is defined by $A \mapsto A(n):= A \otimes \un{\Z} (n)$, 
where $A(n) = \left( A(n)_i, q(n)_{i,j}\right)$ with
$$\tuborg A(n)_i := A_{i+ 2n} \otimes \left( \Spec (k), n\right), \mbox{ and } \\
q(n)_{i,j} = q_{i+2n, j+2n} \otimes 1_{\left( \Spec (k), n \right)}.
\sluttuborg$$
\begin{defn}\label{defn:MTM}
We define the category of \emph{mixed Tate motives} ${\sM}{\sT}{\sM}(k;m)$ over
$k[t]/{(t^{m+1})}$ to be the smallest thick subcategory of ${\sD}{\sM}(k;m)$ 
containing all Tate objects $\un{\Z} (n)$.
\end{defn}
It is clear from our construction and the definition of a thick subcategory
of a monoidal triangulated category (\emph{cf.} \cite[p. 424]{Levine})
that ${\sM}{\sT}{\sM}(k;m)$ is in fact a tensor triangulated category and its
objects are those in ${\sD}{\sM}(k;m)$ which have finite filtrations  
whose graded pieces are the direct sums of Tate objects.

\subsection{Comparison with ${\sD}{\sM}(k)$}
We have seen before (as is obvious from the construction) that 
${\sD}{\sM}(k;0)$ is canonically isomorphic to the integral version of ${\sD}{\sM}(k)$ of Hanamura.
Moreover, if we take ${\sC}'$ to be the partial dg-category 
as before except that we take the morphisms to be only the higher Chow cycles,
i.e., 
\[
\hom_{\mathcal{C}'} ((X,r), (Y,s)) = 
z^{r+s}\left(X \times Y , -\bullet \right),
\] 
then for all $m \geq 1$, there are natural inclusion and forgetful functors
$\iota : {\rm PreTr}(\sC') \to {\rm PreTr}(\sC)$ and
${\rm Forget} :{\rm PreTr}(\sC) \to {\rm PreTr}(\sC')$. These induce the
exact functors
\begin{equation}\label{eqn:Functor*}
\iota : {\sD}{\sM}(k) \to {\sD}{\sM}(k;m), \ \ \mbox{ and }
\end{equation} 
\[
{\rm Forget} : {\sD}{\sM}(k;m) \to 
{\sD}{\sM}(k).
\]
Moreover, for any $X \in \SmProj/k$, there is a split exact sequence
\begin{equation}\label{eqn:exactSeq}
0 \to \TH^r(X,n;m) \to 
\Hom_{{\sD}{\sM}(k;m)}\left(\un{\Z}, h(X)(r)[2r-n]\right) 
\overset{\leftarrow}{\to} \CH^r (X, n) \to 0.
\end{equation}

For a more general smooth quasi-projective variety $X$, not necessarily 
projective, there is a similar split exact sequence, where 
$\TH^r (X, n;m)$ is replaced by the logarithmic additive Chow group 
$\TH^r _{\log} (X, n;m)$ of \cite{KL}, and the functor $h(-)$ is replaced by 
a functor $bm(-)$ that extends $h(-)$ on $\SmProj/k$ to more general 
$k$-schemes, whose construction is the goal of the next section.

\section{Motives of schemes}\label{section:motive-of-schemes}
In this section, we extend the homological functor $h: {\SmProj}/k
\to {\sD}{\sM}(k;m)$ to the category of schemes of finite type over
$k$, assuming the resolution of singularities in the sense of Hironaka.
This will complete the proof of our main theorem.
Using the additive cycle complexes of objects of ${\sD}{\sM}(k;m)$, this extension allows us to get directly an additive cycle complex associated to any scheme $X$
whose homology is the logarithmic additive Chow groups 
$\TH^r_{\log}(X, n ;m)$ of \cite[Theorem~3.3]{KL}. Since our extension
of the homological functor heavily uses the intermediate category 
${\sD}_{\hom}(k)$ of \cite[Section~2]{KL}, we begin this section by recalling
its definition and the related concepts. Throughout this section, we assume 
that the ground field $k$ admits Hironaka's resolution of singularities.  

Let $\Z\SmProj/k$ be the additive category generated by $\SmProj/k$: for any integral smooth projective varieties $X, Y$, define
\[
\Hom_{\Z\SmProj/k}(X,Y):=\Z[\Hom_{\SmProj/k}(X,Y)]
\]
and extend to finite formal sums of integral smooth projective varieties in the natural way. The composition law in $\Z\SmProj/k$ is induced from $\SmProj/k$.

We form the category of bounded complexes $C^b(\Z\SmProj/k)$ and the homotopy 
category $K^b(\Z\SmProj/k)$. We denote the complex concentrated in degree 0 
associated to $X\in\SmProj/k$ by $[X]$. Sending $X$ to $[X]$ defines the 
functor
\[
[-]:\SmProj/k\to C^b(\Z\SmProj/k)
\]

Let $i:Z\to X$ be a closed immersion in $\SmProj/k$, $\mu:X_Z\to X$ the 
blow-up of $X$ along $Z$ and $i_E:E\to X_Z$ the exceptional divisor with the 
structure morphism $q:E\to Z$. Let $C(\mu)$ be the complex
\begin{equation}\label{eqn:BL}
[E]\xrightarrow{(i_E,-q)}[X_Z]\oplus [Z]\xrightarrow{\mu+i}[X]
\end{equation}
with $[X]$ in degree 0. 
\begin{defn}\label{defn:Dhom} The category $\sD_\hom(k)$ is the localization 
of the triangulated category $K^b(\Z\SmProj/k)$ with respect to the thick 
subcategory generated by the complexes $C(\mu)$.

Let 
\[
m_\hom:\SmProj/k\to \sD_\hom(k)
\]
be the composition of functors
\[
\SmProj/k\xrightarrow{[-]}C^b(\SmProj/k)\to K^b(\SmProj/k)\to \sD_\hom(k).
\]
\end{defn}

Recall that $\Sch/k$ is the category of all quasi-projective schemes over $k$ 
and ${\Sch}'/k$ is its subcategory with only proper morphisms. 
\begin{thm}\label{thm:Extension0}{\cite[Theorem~2.9]{KL}} 
The functor $m_\hom$ extends to a functor
\[
M_\hom:\Sch_k'\to  \sD_\hom(k)
\]
such that\\
\\
1. If $\mu:Y\to X$ is a proper morphism in $\Sch_k$, $i:Z\to X$ a closed 
immersion such that $\mu:\mu^{-1}(X\setminus Z)\to X\setminus Z$ is an 
isomorphism, then 
\[
M_\hom(\mu^{-1}(Z)) \to M_\hom(Y) \oplus M_\hom(Z) \to M_\hom(X)
\to M_\hom(\mu^{-1}(Z))[1]
\]
is a distinguished triangle in $\sD_\hom(k)$.
\\
2. Let $j:U\to X$ be an open immersion in $\Sch_k$ with closed complement 
$i:Z\to X$. We have the object $\Cone([i])$ in $C^b(\Z\SmProj/k)$, giving the 
object $m_\hom(\Cone([i]))$ in $\sD_\hom(k)$. Then there is a canonical 
isomorphism
\[
M_\hom(U)\cong m_\hom(\Cone([i]))
\]
in $\sD_\hom(k)$,
giving  a canonical distinguished triangle
\[
M_\hom(Z)\xrightarrow{i_*}M_\hom(X)\xrightarrow{j^*}M_\hom(U)\to M_\hom(Z)[1]
\]
in $\sD_\hom(k)$, natural with respect to proper morphisms of pairs 
$f:(X,U)\to (X',U')$.
\end{thm}

Next we have the following variant of \cite[Proposition~5.5]{Ha1}
for our category ${\sD}{\sM}(k;m)$.
\begin{lem}\label{lem:Yoneda}
Let $u :K \to L$ be a morphism in ${\sD}{\sM}(k;m)$ such that for any $X\in \SmProj/k$ and $s, i \in \Z$, the map
\[
u \circ (-) : \Hom_{{\sD}{\sM}(k;m)}\left((X,s)[i], K\right) \to
\Hom_{{\sD}{\sM}(k;m)}\left((X,s)[i], L\right)  
\]
is an isomorphism. Then $u$ is an isomorphism.
\end{lem}
\begin{proof} This is a straightforward consequence of Yoneda's Lemma
by a spectral sequence argument in \emph{loc. cit.}
\end{proof} 

For $X\in \SmProj/k$ and $r, i \in \Z$, let $h(X)(r)[i]$ 
denote the object $(X,r)[i]$ of ${\sD}{\sM}(k;m)$.
\begin{lem}\label{lem:Blow-up}
Let 
\[
\xymatrix{
E\ar[r]^{i_E}\ar[d]_q&X_Z\ar[d]^\mu\\
Z\ar[r]_i&X
}
\]
be a blow-up square in $\SmProj/k$. Then
\[
h(E)(r)\xrightarrow{(i_{E*},-q_*)}h(X_Z)(r)\oplus 
h(Z)(r)\xrightarrow{\mu_*+i_*}h(X)(r)
\]
is an exact triangle in ${\sD}{\sM}(k;m)$.
\end{lem}
\begin{proof} Let $C(\mu)$ be the complex in ~\eqref{eqn:BL}. It suffices to
show that it is isomorphic to the zero object in ${\sD}{\sM}(k;m)$.
By Lemma~\ref{lem:Yoneda}, it suffices to show that 
$\Hom_{{\sD}{\sM}(k;m)}\left((Y,s)[i], C(\mu) \right)$ is zero for
all $Y \in {\SmProj}/k$. Since $Y \times X_Z$ is the blow-up of $Y \times X$
along $Y \times Z$, it suffices to show that 
$z^r\left(C(\mu), -\bullet ; m \right)$ is acyclic for arbitrary blow-up
$X_Z \to X$ in $\SmProj/k$ and $r \in \Z$. But this follows directly from the 
definition of $z^r\left(C(\mu), -\bullet ; m \right)$, the blow-up formula
for higher Chow groups (\emph{cf.} \cite[Lemma~5.7]{KL}) and the blow-up
formula for the additive higher Chow groups (
\cite[Theorem~5.8]{KL}, \cite[Theorem~3.2]{KP}).
\end{proof}
\begin{prop}\label{prop:extension} The functor $h : \SmProj/k \to 
{\sD}{\sM}(k;m)$ canonically extends to an exact functor
${\sD}(h): \sD_\hom(k) \to {\sD}{\sM}(k;m)$ of triangulated categories.
\end{prop}
\begin{proof}  
Assuming that we can canonically extend the functor $h$ to a functor 
$K^b(h): K^b(\Z\SmProj/k) \to {\sD}{\sM}(k;m)$, the proposition follows from
the description of $\sD_\hom(k)$ in Definition~\ref{defn:Dhom} and
Lemma~\ref{lem:Blow-up}. So we only need to construct $K^b(h)$.

It follows from the definition of the shift functor and the standard 
distinguished triangles in ${\sD}{\sM}(k;m)$ 
(\emph{cf.} proof of Proposition~\ref{prop:TriangC})
if $M(n) : = 
(X_0 \xrightarrow{f_0} X_1 \xrightarrow{f_1} \cdots 
\xrightarrow{f_n} X_{n+1})$ is an object of $K^b(\Z\SmProj/k)$, then
\[
M(n) \cong \Cone\left(M(n-1) \xrightarrow{f_n} X_{n+1}[-n]\right).
\]
Hence, it suffices to show by induction on the length of $M$ that 
$K^b(h) (M)$ is a well defined object of ${\sD}{\sM}(k;m)$.
Now if $M = (X_0 \xrightarrow{f} X_1)$, then $f$ is represented by
a cycle $f \in z^r(X_0 \times X_1, 0)$ which implies that $d(f) = 0$.
Since $h(X_i) = ((X_i, 0), 0)$ for $i = 0,1$, we see that $D(K^b(h)(f)) = 0$
in ${\sD}{\sM}(k;m)$. In particular, the definition of the cone
of a morphism as given above implies that $K^b(M) = (X_0 \xrightarrow{f} X_1)$
is a twisted complex $K$ with $K_i = X_i$ for $i = 0,1$ and $q_{0,1} = f$ 
and hence defines a unique object of ${\sD}{\sM}(k;m)$. Thus $K^b(h)$
canonically extends the functor $h$ from $\SmProj/k$ to an exact functor
$K^b(\Z\SmProj/k) \to {\sD}{\sM}(k;m)$. This completes the proof of the
proposition.
\end{proof}  
\begin{thm}\label{thm:extension1}
The functor $h: \SmProj/k \to {\sD}{\sM}(k;m)$ extends canonically to a
functor
\[
bm : \Sch'/k \to {\sD}{\sM}(k;m)
\]
such that\\
\\
1. If $\mu:Y\to X$ is a proper morphism in $\Sch_k$, $i:Z\to X$ a closed 
immersion such that $\mu:\mu^{-1}(X\setminus Z)\to X\setminus Z$ is an 
isomorphism, then 
\[
bm(\mu^{-1}(Z)) \to bm(Y) \oplus bm(Z) \to bm(X)
\to bm(\mu^{-1}(Z))[1]
\]
is a distinguished triangle in ${\sD}{\sM}(k;m)$.
\\
2. If $j:U\to X$ is an open immersion in $\Sch_k$ with closed complement 
$i:Z\to X$, then there is a canonical distinguished triangle
\[
bm(Z)\xrightarrow{i_*} bm(X)\xrightarrow{j^*}bm(U)\to bm(Z)[1]
\]
in ${\sD}{\sM}(k;m)$, natural with respect to proper morphisms of pairs 
$f:(X,U)\to (X',U')$.
\end{thm}
\begin{proof}
This is an immediate consequence of Theorem~\ref{thm:Extension0} and
Proposition~\ref{prop:extension}.
\end{proof}

\begin{thm}\label{thm:cycle-motive}
The functor $(X, r) \mapsto \sZ \left((X,r), -\bullet ; m\right) =z^r(X, -\bullet ; m)$ canonically extends
to an exact functor $\sZ(-, -\bullet ; m): {\sD}{\sM}(k;m) \to D^{-}(\Z)$ and 
hence defines the mixed cycle complexes of all motives such that for 
$X \in \Sch'/k$, one has 
\[
H^{-i}\left(\sZ\left( (X,r), -\bullet ; m \right)
\right) = {\rm CHC}^r (X, i)\oplus \TH^r_\log(X, i;m),
\]
where ${\rm CHC}^r (X, i) : =  
\Hom_{{\sD}{\sM}(k)}\left(\un{\Z}, bm(X)(r)[2r-n]\right)$
(\emph{cf.} \cite[Definition~2.4]{Ha2}, \cite[Definition~4.4]{Ha1}). 
\end{thm}  

\begin{proof}
Let $K = \left(K^i, q_{i,j}\right)$ be a twisted complex in
${\sD}{\sM}(k;m)$. We define the mixed cycle complex of $K$ following
\cite{Ha1}, using our refined moving lemma for additive cycle complexes.
For each $\alpha \in I(i)$, we take a distinguished subcomplex 
$\sZ(K^i_{\alpha},- \bullet ; m)' \subset \sZ(K^i_{\alpha}, -\bullet ; m)$
so that letting $\sZ(K^i, -\bullet ; m)' := 
\stackrel{}{\underset {\alpha \in I(i)}{\oplus}} 
\sZ(K^i_{\alpha}, -\bullet ; m)'$, the map 
\[
\left(q_{i_{r-1}, i_r} \circ \cdots \circ q_{i_{r_0}, i_{r_1}}\right)_*
: \sZ(K^i, -\bullet ; m)' \to \sZ(K^j, -\bullet +j-i-r ; m)'
\]
is defined and associative for any sequence $i = i_0 < \cdots < i_r = j$.  
The \emph{mixed cycle complex} $\sZ(K, -\bullet ; m)$ of $K$ is defined as the 
complex $(L, d)$ with
\[
L^i = \stackrel{}{\underset {j}{\bigoplus}}
\stackrel{}{\underset {\alpha \in I(j)}{\bigoplus}} 
z^{r_{\alpha}}(K^j_{\alpha}, j-i ; m)',
\]
\[
d^i =  \stackrel{}{\underset {j}{\sum}}
\left((-1)^j \delta_j + \stackrel{}{\underset {j < l}{\sum}}
{(q_{j,l})}_ * \right).
\]
This defines the desired functor $\sZ\left(-, -\bullet ; m\right)$.

It is easy to check from this definition that if $M =
(X_0 \xrightarrow{f_0} X_1 \xrightarrow{f_1} \cdots 
\xrightarrow{f_n} X_{n+1})$ is an object of $K^b(\Z\SmProj/k)$, then
$\sZ\left(M(r), -\bullet ; m\right)$ is the total complex 
associated to the double complex given by 
\[
{\sZ\left(M(r), -\bullet ; m\right)}^{i,j}
= z^r(X_i, j; m)'
\]  
with the horizontal differential given by $(-1)^i{(f_i)}_*$ and the vertical
differential given by ${\delta}_{X_j}$.

If $X$ is now a complete variety (possibly singular) with a smooth cubical
resolution $X_{\bullet} \to X$, let $M = X_{\bullet}$ also denote the 
associated chain complex in $K^b(\Z\SmProj/k)$ with the differential being the 
alternating sum of the face maps of the cubical object $X_{\bullet}$.
So, the above description of $\sZ\left(M(r), -\bullet ; m\right)$ and
\cite[Theorem~6.1]{KL} immediately imply that 
$H^{-i}\left(\sZ\left(bm(X)(r), -\bullet ; m \right)\right) = 
%H^i\left(\sZ\left(X_{\bullet}(r), i ; m \right)\right) =
\CH^r (X, i) \oplus \TH^r_\log(X, i;m)$. If $X$ is not complete, the corresponding
isomorphism now follows from Theorem~\ref{thm:extension1}
and \cite[Corollary~6.2]{KL}. This completes the proof of the theorem. 
\end{proof}

\begin{defn}{(Total higher Chow groups of a motive)}\label{defn:CCM}
For $A \in {\sD}{\sM}(k;m)$, we define its \emph{motivic 
cohomology} by 
\begin{equation}\label{eqn:CCM1}
\CH_{\rm log}(A,n;m): = H^n \left(\sZ(A, -\bullet;m)\right),
\end{equation}
where $\sZ(-, \bullet ;m)$ is the functor of Theorem~\ref{thm:cycle-motive}.
\end{defn}
\begin{cor}\label{cor:CCM2}
For a smooth quasi-projective variety $X$, one has 
\[
\Hom_{{\sD}{\sM}(k;m)}\left(\un{\Z}, bm(X)(r)[2r-n]\right) =
\CH^r(X,n) \oplus \TH^r_{\rm log}(X,n;m).
\]
\end{cor}
\begin{proof} For $X$ smooth and projective, this is shown below. If $X$
is a projective but possibly singular variety, let $X_{\bullet} \to X$ be
a smooth cubical resolution. Since 
$\Hom_{{\sD}{\sM}(k;m)}\left(\un{\Z}, -\right)$ is a cohomological functor,
the spectral sequence
\[
E^{p,q}_1 = \Hom_{{\sD}{\sM}(k;m)}\left(\un{\Z}, h(X_p)[p-q]\right)
\Rightarrow \Hom_{{\sD}{\sM}(k;m)}\left(\un{\Z}, bm(X)\right)
\]
and Theorem~\ref{thm:cycle-motive} show that 
\begin{equation}\label{eqn:CCM3}
\Hom_{{\sD}{\sM}(k;m)}\left(\un{\Z}, bm(X)(r)[2r-n]\right) =
{\rm CHC}^r(X,n) \oplus \TH^r_{\rm log}(X,n;m).
\end{equation}
If $X$ is not necessarily projective, ~\eqref{eqn:CCM3} now follows from
Theorem~\ref{thm:extension1}(2). Finally, for $X$ smooth and quasi-projective,
one has ${\rm CHC}^r(X, n) = \CH^r(X,n)$ by \cite[p.328]{Ha2}.
\end{proof}
\noindent
{{\bf{Proof of Theorem~\ref{thm:Main}:}}} The part $(1)$ of the theorem is 
already shown in ~\eqref{eqn:Functor*} and ~\eqref{eqn:exactSeq}.

For $(2)$, we first observe from the definition of the differential of
a complex $\left(\Hom_{{\rm PreTr}(\sC)}(A, B), D\right)$ in 
~\eqref{eqn:D} that
for twisted complexes $A,B$ with $p_{i,j}$, $q_{i,j}$'s all zero, one has $D = d$.  
Hence, we see from ~\eqref{eqn:mfunctor} and ~\eqref{subsection:Ttwist} that
 \[
\Hom_{{\sD}{\sM}(k;m)}\left(\un{\Z}, h(X)(r)[2r-n]\right)
= H^n\left({\sZ}^r(X, -\bullet ;m)\right).
\]
Part $(2)$ now follows from ~\eqref{eqn:Pdcom*}. The last part is shown
in Corollary~\ref{cor:CCM2}.
\\
\\
\noindent\emph{Acknowledgments.} The second author would like to thank TIFR 
and KAIST for their hospitality and administrative help. He wishes to thank 
Juya for her moral support. During this work, JP was partially supported by 
Basic Science Research Program through the National Research Foundation of 
Korea (NRF) funded by the Ministry of Education, Science and Technology 
(2009-0063180).


\begin{thebibliography}{99}

\bibitem{BS} Balmer, P., Schlichting, M., {\sl Idempotent completion of 
triangulated categories}, J. Algebra, \textbf{236,} (2001), no. 2, 819-834.

\bibitem{Bloch-tangent} Bloch, S., {\sl On the tangent space to Quillen 
$K$-theory\/}, Algebraic $K$-theory, I: Higher $K$-theories 
(Proc. Conf., Battelle Memorial Inst., Seattle, Wash., 1972), pp. 205-210. 
Lecture Notes in Math., \textbf{341,} Springer, Berlin, 1973.

\bibitem{Bl1} Bloch, S., {\sl Algebraic cycles and higher $K$-theory\/},
Adv. Math., \textbf{61,} (1986), 267-304.

\bibitem{Bloch} Bloch, S., {\sl Algebraic cycles and the Lie algebra of mixed 
Tate motives\/},
J. Amer. Math. Soc., \textbf{4,} (1991), no. 4, 771-791.

\bibitem{Bl2} Bloch, S., {\sl The moving lemma for higher Chow groups\/}, 
J. Algebraic Geom., \textbf{3,} (1994), no. 3, 537-568.

\bibitem{BE2} Bloch, S., Esnault, H., {\sl The additive dilogarithm\/},
Doc. Math., J., DMV \textbf{Extra Vol.,} (2003), 131-155.

\bibitem{BK} Bondal, A. I.; Kapranov, M. M., {\sl Framed triangulated 
categories}, (Russian) Mat. Sb. \textbf{181,} (1990), no. 5, 669-683; 
translation in English as {\sl Enhanced triangulated categories}, Math. 
USSR-Sb. \textbf{70,}(1991), no. 1, 93-107.

\bibitem {Fulton} Fulton, W., {\sl Intersection theory\/},
Second Edition, Ergebnisse der Mathematik und ihrer Grenzgebiete 3,
Folge. A Series of Modern Surveys in Mathematics, {\bf 2,}
Springer-Verlag, Berlin, 1998.

\bibitem{Goncharov} Goncharov, A., {\sl Euclidean scissor congruence groups 
and mixed Tate motives over dual numbers}, Math. Res. Lett., \textbf{11,} 
(2004), no. 5-6, 771-784.

\bibitem{GN} Guill\'en, F., Navarro Aznar, V., 
{\sl Un crit\`ere d'extension des foncteurs 
d\'efinis sur les sch\'emas lisses\/},
Publ. Math. Inst. Hautes \'Etudes Sci., {\bf 95,} (2002), 1-91.

\bibitem {Ha2} Hanamura, M., {\sl Homological and cohomological motives of
algebraic varieties}, Invent. Math., \textbf{142,} (2000) 319-149.

\bibitem {Ha1} Hanamura, M., {\sl Mixed motives and algebraic cycles II}, 
Invent. Math., \textbf{158,} (2004) 105-179.

\bibitem{Hesselholt} Hesselholt L., {\sl $K$-theory of truncated polynomial 
algebras\/}, Handbook of K-Theory, {\bf 1,}
Springer-Verlag, Berlin, 2005,  71-110.

\bibitem{KL} Krishna, A., Levine, M., {\sl Additive higher Chow groups of 
schemes\/}, J. Reine Angew. Math., \textbf{619,} (2008) 75-140.

\bibitem{KP} Krishna, A., Park, J., {\sl Moving lemma for additive Chow 
groups and applications}, Preprint Sep. 2009, arXiv:0909.3155

\bibitem{KP1} Krishna, A., Park, J., {\sl Additive higher Chow group of 
1-cycles on fields\/}, In preparation. 
 
\bibitem{Levine} Levine, M., {\sl Mixed Motives\/},
Mathematical Surveys and Monographs, {\bf 57},
American Mathematical Society, Providence, RI, 1998.
 

\bibitem{P1} Park, J., {\sl Regulators on additive higher Chow groups\/},
Amer. J. Math., \textbf{131,} (2009) no. 1, 257-276.

\bibitem{P2} Park, J., {\sl Algebraic cycles and additive dilogarithm\/},
Int. Math. Res. Not., \textbf{2007,} no. 18, Article ID rnm067.

\bibitem{R} R\"ulling, K., {\sl The generalized de Rham-Witt complex over a 
field is a complex of zero-cycles\/}, 
J. Algebraic Geom., \textbf{16,} (2007), no. 1, 109-169.

\bibitem{Soule} Soul\'e, C., 
{\sl Op\'erations en $K$-th\'eorie alg\'ebrique\/},
Canad. J. Math.,  \textbf{37,}  (1985),  no. 3, 488-550. 

\bibitem{Voevodsky} Voevodsky, V., {\sl Triangulated categories of motives 
over a field.} in Cycles, Trans- fers, and Motivic Cohomology Theories, ed. 
by V. Voevodsky, A. Suslin, E.M. Friedlander. Princeton: Princeton Univ. 
Press 2000
\end{thebibliography}
\end{document}